\documentclass[12pt,reqno]{amsart}
\theoremstyle{plain}
\newtheorem{thm}{Theorem}[section]
\numberwithin{equation}{section} %% Comment out for sequentially-numbered

\newtheorem{cor}[thm]{Corollary} %%Delete [thm] to re-start numbering
\newtheorem{lemma}[thm]{Lemma} %%Delete [thm] to re-start numbering
\newtheorem{prop}[thm]{Proposition}
\newtheorem{defi}[thm]{Definition}
\newtheorem{remark}[thm]{Remark}
\newtheorem{examples}[thm]{Examples}
\theoremstyle{remark}

\numberwithin{figure}{section} %% Comment out for sequentially-numbered

\usepackage{amscd,amssymb,euscript}
\newcount\theTime
\newcount\theHour
\newcount\theMinute
\newcount\theMinuteTens
\newcount\theScratch
\theTime=\number\time \theHour=\theTime \divide\theHour by 60
\theScratch=\theHour \multiply\theScratch by 60
\theMinute=\theTime \advance\theMinute by -\theScratch
\theMinuteTens=\theMinute \divide\theMinuteTens by 10
\theScratch=\theMinuteTens \multiply\theScratch by 10
\advance\theMinute by -\theScratch

\def\today{{\number\day\space
 \ifcase\month\or
  January\or February\or March\or April\or May\or June\or
  July\or August\or September\or October\or November\or December\fi
 \space\number\year}}

\newcommand\Bc{{\mathcal{B}}}

\newcommand\Cpx{{\mathbf C}}

\newcommand\DEu{{\EuScript D}}

\newcommand\dif{\mbox{\it d}}

\newcommand\Ec{{\mathcal{E}}}

\newcommand\eps{\epsilon}

\newcommand\Fc{{\mathcal{F}}}

\newcommand\fs{_{\text{\it fs}}}

\newcommand\Gc{{\mathcal{G}}}

\newcommand\HEu{{\EuScript H}}                   % requires package euscript

\newcommand\htil{{\tilde h}}

\newcommand\Ic{{\mathcal{I}}}

\newcommand\ImagPart{{\mathrm{Im}\;}}

\newcommand\Ints{{\mathbf Z}}

\newcommand\Jc{{\mathcal{J}}}

\newcommand\Kc{{\mathcal{K}}}

\newcommand\Lc{{\mathcal{L}}}

\newcommand\Mcal{{\mathcal{M}}} %  needed to be renamed, because was ``\Mc already defined''
\newcommand\Mcalb{{\overline\Mcal}}

\newcommand\Nats{{\mathbf N}}

\newcommand\Qt{{\widetilde Q}}

\newcommand\RealPart{{\mathrm{Re}\;}}

\newcommand\Reals{{\mathbf R}}

\newcommand\subb{_{\text{\it b}}}

\newcommand\Tt{{\widetilde T}}

\def\supp{\text{supp }}

\begin{document}

\pagestyle{myheadings}

\title{Sums of commutators in ideals and modules of type II factors}

\author{K.J.\ Dykema}

\address{\hskip-\parindent
K.J.\ Dykema\\
Department of Mathematics\\
Texas A\&M University\\
College Station TX 77843--3368, USA}
\email{Ken.Dykema@math.tamu.edu}

\author{N.J.\ Kalton}

\address{\hskip-\parindent
N.J.\ Kalton \\
Department of Mathematics \\
University of Missouri \\
Columbia MO  65211, USA} \email  {\tt nigel@math.missouri.edu}

\thanks{K.J.D.\ supported in part by NSF grant DMS--0070558.
N.J.K.\ supported in part by NSF grant DMS--9870027.}

\begin{abstract}
Let $\Mcal$ be a factor of type II$_\infty$ or II$_1$ having
separable predual and let $\Mcalb$ be the algebra of affiliated
$\tau$--measureable operators. We characterize the commutator
space $[\Ic,\Jc]$ for sub--$(\Mcal,\Mcal)$--bimodules $\Ic$ and
$\Jc$ of $\Mcalb$.
\end{abstract}

\date{21 April, 2003}

\maketitle

%\markboth{\tiny Type II commutators \timeanddate}{\tiny Type II commutators \timeanddate}

\section{Introduction and description of results}

Let $\Mcal$ be a von Neumann algebra of type II$_\infty$ having
separable predual. We will study the commutator structure of
ideals of $\Mcal$ and, more generally, of modules of operators
affiliated to $\Mcal$.

Fix a faithful semifinite trace $\tau$ on $\Mcal$, and let $\Mcal$
be represented on a Hilbert space $\HEu$. Segal~\cite{Segal}
introduced measurability for unbounded operators on $\HEu$
affiliated to $\Mcal$. Later Nelson~\cite{Nelson}, in a slightly
different approach, defined the completion $\Mcalb$ of $\Mcal$
with respect to a notion of convergence in measure, and showed
that the operations on $\Mcal$ extend to make $\Mcalb$ a
topological $*$--algebra. He also showed that $\Mcalb$ is the set
of all {\em $\tau$--measurable} operators, i.e. the closed,
densely defined, possibly unbounded operators $T$ on $\HEu$,
affiliated with $\Mcal$, such that for every $\eps>0$ there is a
projection $E\in\Mcal$ with $\tau(1-E)<\eps$ and with $TE$
bounded. Note that $\Mcalb$ is defined independently of the
Hilbert space $\HEu$ on which $\Mcal$ acts, but is then
characterized in terms of operators on $\HEu$. Nelson's work was
done in the more general context of a von Neumann algebra $\Mcal$
equipped with a fixed finite or semifinite faithful normal trace.
(See~\cite{Cecchini} for a proof that Segal's and Nelson's
definitions are equivalent in II$_\infty$ factors.)

We consider subspaces $\Ic\subseteq\Mcalb$ that are globally
invariant under left and right multiplication by elements from
$\Mcal$. These are thus sub--($\Mcal,\Mcal$)--bimodules of
$\Mcalb$; for brevity we will call them {\em submodules} of
$\Mcalb$. Note that if such a submodule $\Ic$ is actually
contained in $\Mcal$, then it is a two--sided ideal of $\Mcal$.
Submodules of $\Mcalb$ are analogues in the type II$_\infty$
context of ideals of $B(\HEu)$ in the type I context. The
submodules of $\Mcalb$ can be classified in terms of the singular
numbers of their elements, analogously to Calkin's
classification~\cite{Calkin} of the ideals of $B(\HEu)$. If
$T\in\Mcalb$ and $t>0$, the \emph{$t$-th singular number} of $T$
is
\begin{equation}\label{eq:mut}
\mu_t(T)=\inf\big(\{\|T(1-E)\|:\,E\in\Mcal \text{ a projection
with }\tau(E)\le t\}\cup\{0\}\big)\;,
\end{equation}
and we denote by $\mu(T)$ the function $t\mapsto\mu_t(T)$. If
$\Ic\subseteq\Mcalb$ is a submodule, we set
\[
\mu(\Ic)=\{\mu(T)\mid T\in\Ic\}
\]
and we call $\mu(\Ic)$ the \emph{characteristic set} of $\Ic$. The
aforementioned classification is the bijection
$\Ic\mapsto\mu(\Ic)$ from the set of all submodules of $\Mcalb$ to
the set of all characteristic sets, where, abstractly, a
characteristic set is a set of decreasing functions on
$(0,\infty)$ satisfying certain properties. Several authors have
used singular numbers to characterize ideals of $\Mcal$ and
modules of $\Mcalb$ (see~\cite{DoddsDoddsPagter},
\cite{StrohWest}, \cite{West}
and~\cite{DoddsPagterSemenevSukochev}), and the full
classification result was derived by Guido and Isola
in~\cite{GuidoIsola}.

One interesting facet of submodules of $\Mcalb$ is that their
classification involves both asymptotics at infinity (the rate of
decay of $\mu_t(T)$ as $t\to\infty$) and asymptotics at zero (the
rate of increase of $\mu_t(T)$ as $t\to0$).

We consider additive commutators $[A,B]=AB-BA$ of elements of
$\Mcalb$ and study the commutator spaces
\[
[\Ic,\Jc]=\{\sum_{k=1}^n[A_k,B_k]\mid
n\in\Nats,\,A_k\in\Ic,\,B_k\in\Jc\}
\]
of submodules $\Ic$ and $\Jc$ of $\Mcalb$. Note
$[\Ic,\Jc]\subseteq\Ic\Jc$, where $\Ic\Jc$ is the submodule of
$\Mcalb$ spanned by all products $AB$ with $A\in\Ic$ and
$B\in\Jc$. Using properties of singular numbers (which are
reviewed in~\S2), one easily shows that $\mu(\Ic\Jc)$ is the set
of all decreasing functions $f:(0,\infty)\to[0,\infty)$ bounded
above by products $gh$ with $g\in\mu(\Ic)$ and $h\in\mu(\Jc)$.
Since an element of $\Ic\Jc$ belongs to $[\Ic,\Jc]$ if and only if
its real and imaginary parts belong to $[\Ic,\Jc]$, to characterize
$[\Ic,\Jc]$ it will suffice to describe the normal elements of it.
This we do as follows: given a normal element $T\in\Ic\Jc$, let
$E_{|T|}$ denote the spectral measure of the positive part $|T|$
of $T$. Then $T\in[\Ic,\Jc]$ if and only if there is
$h\in\mu(\Ic\Jc)$ such that
\begin{equation}\label{eq:Tcomm}
|\tau(TE_{|T|}(\mu_s(T),\mu_r(T)])|\le rh(r)+sh(s)
\end{equation}
for all $0<r<s<\infty$. This is analogous, though for asymptotics
in both directions, to the characterization of commutator spaces
for ideals of $B(\HEu)$ found in~\cite{DFWW} (see
also~\cite{Kalton-L1} for the earlier result in the case of the
trace--class operators). Our proof relies on a result of Fack and
de la Harpe~\cite{FackdelaHarpe}, expressing any trace--zero
element of a II$_1$--factor as a sum of a fixed number of
commutators of elements whose norms are controlled. A corollary of
our characterization is
\[
[\Ic,\Jc]=[\Ic\Jc,\Mcal]
\]
for any submodules $\Ic$ and $\Jc$ of $\Mcalb$. We also give a
characterization of $T\in[\Ic,\Jc]$ for $T$ normal that considers
separately the asymptotics at $0$ and at $\infty$.

As an alternative to using the characteristic set $\mu(\Ic)$ of a
submodule $\Ic\subseteq\Mcalb$ for the classification of
submodules, one can use the corresponding rearrangement invariant
function space $S(\Ic)$, which is the set of all measureable
functions $f:(0,\infty)\to\Cpx$ such that the decreasing
rearrangement of the absolute value of $f$ lies in $\mu(\Ic)$.
Then every normal element $T\in\Ic$ gives rise to a unique (up to
rearrangement) function $f_T\in S(\Ic)$ defined as follows: Fix
any measure preserving transformation from $(0,\infty)$ with
Lebesgue measure to the disjoint union of four copies of
$(0,\infty)$ with Lebesgue measure, in order to define the
measureable function $g_1\oplus g_2\oplus g_3\oplus
g_4:(0,\infty)\to\Cpx$, given measureable functions
$g_j:(0,\infty)\to\Cpx$. Now let
$f_T=f_1\oplus(-f_2)\oplus(if_3)\oplus(-if_4)\in S(\Ic)$, where
\begin{alignat*}{2}
f_1(t)&=\mu_t((\RealPart T)_+)\quad&f_2(t)&=\mu_t((\RealPart T)_-) \\
f_3(t)&=\mu_t((\ImagPart T)_+)\quad&f_4(t)&=\mu_t((\ImagPart
T)_-)\;,
\end{alignat*}
with $\RealPart T=(T+T^*)/2=(\RealPart T)_+-(\RealPart T)_-$,
where $(\RealPart T)_+$ and $(\RealPart T)_-$ are commuting
positive operators whose product is zero, and similarly for
$\ImagPart T=(T-T^*)/2i=(\ImagPart T)_+-(\ImagPart T)_-$.
Then in the case when $\lim_{t\to\infty}\mu_t(T)=0$ for all elements $T\in\Ic\Jc$, the
condition~\eqref{eq:Tcomm} above for $T\in[\Ic,\Jc]$ with $T$
normal can be rephrased in terms of $f_T$ and is seen to be
equivalent to the condition found in~\cite{FigielKalton} for $f_T$
to belong to the kernel of every symmetric functional on $S(\Ic\Jc)$.
Thus, our main result can be seen as a noncommutative analogue of
this result from~\cite{FigielKalton}.
See also~\cite{DoddsPagterSemenevSukochev} for related results on
Banach symmetric functions spaces and the corresponding submodules of $\Mcalb$.

In the case of a II$_1$--factor $\Mcal$, we give an analogous
characterization of the commutator spaces $[\Ic,\Jc]$ for
submodules $\Ic$ and $\Jc$ of $\Mcalb$.

In the case of ideals in $\mathcal B(\mathcal H)$ it was shown in
\cite{Kreine} that for  quasi-Banach ideals $\mathcal I$ the
subspace $[\mathcal I,\mathcal B(\mathcal H)]$ can be
characterized purely in spectral terms (see also \cite{Kalton-L1}
for an earlier result in this direction). More generally this
result was established for the class of geometrically stable
ideals. This means that for such ideals if two operators $S,T$ in
$\mathcal I$ have the same spectrum (counting algebraic
multiplicities) and $S\in [\mathcal I,\mathcal B(\mathcal H)]$
then $T\in [\mathcal I,\mathcal B(\mathcal H)].$ This was known
for hermitian operators (and hence normal operators) from the
results in \cite{DFWW}, but is generally false (see \cite
{DykemaKalton}).  We study the same phenomenon in type
$II_{\infty}-$factors.   In this case, since we need a notion
corresponding to multiplicity we employ the Brown measure
\cite{Brown} as a substitute for the notion of spectrum. The Brown
measure of an operator is a measure with support contained in its spectrum. It
is, however, only defined for certain special types of operators.
Nevertheless we obtain a quite satisfactory analogue of the result
of \cite{Kreine}.  If $\mathcal I$ is a geometrically stable
submodule of $\overline{\mathcal M}$ and $T\in\mathcal I$ admits a
Brown measure $\nu_T$ then $T\in [\mathcal I,\mathcal M]$ if and
only if there is a positive operator $V\in\mathcal I$ so that
$$ \bigg|\int_{r<|z|\le s}z\,d\nu_T(z)\bigg|\le
r\tau(E_{V}(r,\infty))+s\tau(E_V(s,\infty)) \qquad 0<r<s<\infty.$$
This condition depends  only on the Brown measure associated to
$T$.

The paper is organized as follows: In~\S\ref{sec:s-numb}, we
recall some facts about singular numbers of elements of $\Mcalb$.
In~\S\ref{sec:class}, we describe the classification of submodules
of $\Mcalb$ when $\Mcal$ is a type II$_\infty$ or II$_1$ factor
with separable predual. In~\S\ref{sec:comm}, we prove the main
results characterizing $[\Ic,\Jc]$. In~\S\ref{sec:asymp}, we give
a characterization of $[\Ic,\Jc]$ in the II$_\infty$ case,
separating the asymptotics at $0$ and $\infty$.  Results on the
Brown measure are discussed in~\S\ref{spectral}.

\section{Preliminaries on singular numbers}
\label{sec:s-numb}

If $\Mcal$ is a von Neumann algebra with a fixed finite or
semifinite normal trace $\tau$, then the singular numbers
(sometimes called generalized singular numbers) of elements of
$\Mcal$ and more generally of $\tau$--measureable operators
affiliated to $\Mcal$ have been understood for many years; see,
for example,~\cite{MvN}, \cite{Grothendieck-Bourbaki}, \cite{Fack}
and~\cite{FackKosaki}. In this section, we review these concepts
and some results, introduce the notation we will use throughout
the paper and prove a technical result that will be of use later.

Recall that $t$-th singular number of $T\in\Mcalb$ is defined for
$t>0$ by~\eqref{eq:mut}. Since $T$ is $\tau$--measurable, we have
$0\le\mu_t(T)<+\infty$. We will also use the convention
$\mu_0(T)=\|T\|$, where $\|T\|=\infty$ if $T\not\in\Mcal$. Note
that $t\mapsto\mu_t(T)$ is a nonincreasing function from
$[0,\infty)$ into $[0,\infty]$. If $\tau$ is a finite trace, then
by our convention that $\tau(1)=1$ we have $\mu_t(T)=0$ whenever
$t\ge1$. We will use the following properties of singular numbers;
see~\cite{Fack} or~\cite{FackKosaki} for proofs.
\begin{prop}\label{prop:snumb}
Let $\Mcal$ be a von Neumann algebra with a distinguished finite
or semifinite normal faithful trace, let $S,T\in\Mcalb$ and
$s,t\ge0$. Then
\renewcommand{\labelenumi}{(\roman{enumi})}
\begin{enumerate}
\item $\mu_t(T)=\mu_t(T^*)=\mu_t(|T|)$,
\item $\mu_{s+t}(S+T)\le\mu_s(S)+\mu_t(T)$,
\item $\mu_{s+t}(ST)\le\mu_s(S)\mu_t(T)$,
\item if $A,B\in\Mcal$, then $\mu_t(ATB)\le\|A\|\|B\|\mu_t(T)$.
\end{enumerate}
Moreover,
\renewcommand{\labelenumi}{(\roman{enumi})}
\begin{enumerate}
\setcounter{enumi}{3}
\item the function $[0,\infty)\ni t\mapsto\mu_t(T)\in[0,\infty]$ is continuous from the right.
\end{enumerate}
\end{prop}

Given $T\in\Mcalb$, let $A\mapsto E_{|T|}(A)$ be the
projection--valued spectral measure of the positive part $|T|$ of
$T$.
(To avoid clutter, when $A$ is an inverval we will frequently omit to write parenthesis,
writing just $E_{|T|}A$.)
\begin{prop}[\cite{FackKosaki}, 2.2]\label{prop:FK}
For $t\ge0$ we have
\begin{equation}
\label{eq:infatt}
\mu_t(T)=\inf\big(\{s\ge0\mid\tau\big(E_{|T|}(s,\infty)\big)\le
t\}\cup\{\infty\}\big)
\end{equation}
and the infimum is attained, giving
\begin{equation}\label{eq:tauE}
\tau\big(E_{|T|}(\mu_t(T),\infty)\big)\le t
\end{equation}
whenever $\mu_t(T)<\infty$.
\end{prop}

\begin{lemma}\label{lem:infsup}
Let $\Mcal$ be a nonatomic von Neumann algebra with a normal
faithful semifinite trace $\tau$, let $T\in\Mcalb$ and let
$x\in\Reals$, $x\ge0$. Then
\begin{align}
\tau\big(E_{|T|}(x,\infty)\big)&=\inf\big(\{s\ge0\mid\mu_s(T)\le x\}\cup\{\infty\}\big), \label{eq:inf} \\
\tau\big(E_{|T|}[x,\infty)\big)&=\sup\big(\{s\ge0\mid\mu_s(T)\ge
x\}\cup\{0\}\big), \label{eq:sup}
\end{align}
and the infimum in~\eqref{eq:inf} is attained.
\end{lemma}
\begin{proof}
The infimum in~\eqref{eq:inf} is attained because
$s\mapsto\mu_s(T)$ is continuous from the right. If
$a=\tau(E_{|T|}(x,\infty))<\infty$, then, since
$\|T(1-E_{|T|}(x,\infty))\|\le x$, we have $\mu_a(T)\le x$,
proving $\ge$ in~\eqref{eq:inf}. On the other hand, if
$\mu_s(T)\le x<\infty$, then using~\eqref{eq:tauE} we have
\[
\tau\big(E_{|T|}(x,\infty)\big)\le\tau\big(E_{|T|}(\mu_s(T),\infty)\big)\le
s,
\]
proving $\le$ in~\eqref{eq:inf}.

If $s<\tau(E_{|T|}[x,\infty))$, then for any projection
$P\in\Mcal$ with $\tau(P)=s$, we have $(1-P)\wedge
E_{|T|}[x,\infty)\ne0$. Hence $\|T(1-P)\|\ge x$. Therefore
$\mu_s(T)\ge x$, which proves $\le$ in~\eqref{eq:sup}. If
$\tau(E_{|T|}[x,\infty))<s'<\infty$, then since
$[x,\infty)=\bigcap_{0<r<x}(r,\infty)$, there is $r<x$ such that
$\tau(E_{|T|}(r,\infty))\le s'$. But then $\mu_{s'}(T)\le r<x$,
which implies $s'\ge\sup(\{s\ge0\mid\mu_s(T)\ge x\}\cup\{0\})$.
This proves $\ge$ in~\eqref{eq:sup}.
\end{proof}

\begin{defi}\label{defi:oplus}\rm
Let $\Mcal$ be a II$_\infty$--factor and let us introduce the
natural notation $\oplus$. Since $\Mcalb$ consists of (in general
unbounded) operators on a Hilbert space $\HEu$, by choosing an
isomorphism $\HEu\cong\HEu\oplus\HEu$, we may realize
$\Mcalb\oplus\Mcalb$ as a subalgebra of $\Mcalb$ in such a way
that $\tau(S\oplus T)=\tau(S)+\tau(T)$ whenever $S$ and $T$ are in
$L^1(\Mcal,\tau)\subseteq\Mcalb$. Thus for $S,T\in\Mcalb$,
$S\oplus T$ defines an element of $\Mcalb$ uniquely up to
conjugation by a unitary in $\Mcal$. Since $U^*AU=A+[AU,U^*]$
whenever $U$ is unitary and $A\in\Mcalb$, if $S,T\in\Ic$ for any
submodule $\Ic\subseteq\Mcalb$, the direct sum $S\oplus T$ is
defined uniquely up to addition of a commutator from
$[\Ic,\Mcal]$. Moreover, we have $\Ic\oplus\Ic\subseteq\Ic$ and
for every $T\in\Ic$ we get
\[
T\oplus0\in T+[\Ic,\Mcal]
\]
by using an appropriate nonunitary isometry in $\Mcal$.
\end{defi}

\begin{prop}\label{prop:SoplusT}
Let $S,T\in\Mcalb$ and let $a\ge0$. Then
\[
\mu_a(S\oplus T)=\inf\{\max(\mu_b(S),\mu_c(T))\mid
b,c\ge0,\,b+c=a\}\;.
\]
\end{prop}
\begin{proof}
The case $a=0$ is straightforward, so we may assume $a>0$. It is
clearly equivalent to show
\begin{equation}\label{eq:muaST}
\mu_a(S\oplus T)=\inf\{\max(\mu_b(S),\mu_c(T))\mid
b,c\ge0,\,b+c\le a\}\;.
\end{equation}
Given $b,c\ge0$ such that $b+c\le a$, by~\eqref{eq:tauE} we have
\[
\tau(E_{|S|}(\mu_b(S),\infty)\oplus
E_{|T|}(\mu_c(T),\infty))\le b+c\le a\;,
\]
so using the definition~\eqref{eq:mut} of singular numbers, we get
\[
\mu_a(S\oplus T)\le\|SE_{|S|}[0,\mu_b(S)]\oplus
TE_{|T|}[0,\mu_c(T)]\| \le\max(\mu_b(S),\mu_c(T)\;.
\]
This shows $\le$ in~\eqref{eq:muaST}. For the reverse inclusion,
by~\eqref{eq:infatt} we have
\[
\mu_a(S\oplus T)=\inf\{r\ge0\mid\tau(E_{|S\oplus
T|}(r,\infty))\le a\}\;.
\]
But
\[
E_{|S\oplus
T|}(r,\infty)=E_{|S|\oplus|T|}(r,\infty)=E_{|S|}(r,\infty)\oplus
E_{|T|}(r,\infty)\;,
\]
so
\[
\mu_a(S\oplus
T)=\inf\{r\ge0\mid\tau(E_{|S|}(r,\infty))+\tau(E_{|T|}(r,\infty))\le
a\}\;.
\]
By Lemma~\ref{lem:infsup}, if $b=\tau(E_{|S|}(r,\infty))$ and
$c=\tau(E_{|T|}(r,\infty))$, then $\mu_b(S)\le r$ and
$\mu_c(T)\le r$. This implies $\ge$ in~\eqref{eq:muaST}.
\end{proof}

The next lemma can be described as mashing the atoms of $E_{|T|}$.
It is both straightforward and similar to~\cite[Lemma
1.8]{GuidoIsola}. However, for completeness, we include a proof.
\begin{lemma}
\label{lem:mashatoms} Let $\Mcal$ be a von Neumann algebra without
minimal projections and with a distinguished semifinite normal
faithful trace $\tau$. Let $T\in\Mcalb$. Then there is a family
$(P_t)_{t\ge0}$ of projections in $\Mcal$ such that for all $s$
and $t$,
\renewcommand{\labelenumi}{(\roman{enumi})}
\begin{enumerate}
\item $s\le t$ implies $P_s\le P_t$,
\item $\tau(P_t)=t$,
\item $P_t$ and $|T|$ commute, and if $T$ is normal then $P_t$ and $T$ commute,
\item $E_{|T|}(\mu_t(T),\infty)\le P_t\le E_{|T|}[\mu_t(T),\infty)$,
\item if $x>0$, then $E_{|T|}(x,\infty)=P_y$, where $y=\inf\{t>0\mid\mu_t(T)\le x\}$.
\end{enumerate}
Furthermore, suppose $\lim_{t\to\infty}\mu_t(T)=0$. Then, letting
$F$ be the projection--valued Borel measure in $\Mcal$ supported
on $(0,\infty)$ and satisfying $F((a,b))=P_b-P_a$, we have
\[
|T|=\int_{(0,\infty)}\mu_t(T)\dif F(t).
\]
\end{lemma}
\begin{proof}
If $T$ is not normal, then we may replace $T$ by $|T|$, so assume
$T$ is normal. Set $P_t=E_{|T|}(\mu_t(T),\infty)$ whenever
$E_{|T|}(\{\mu_t(T)\})=0$. For these values of $t$, it follows
from Lemma~\ref{lem:infsup} that $\tau(P_t)=t$. The set
\[
\Ec=\{E_{|T|}(\{\mu_t(T)\})\mid t>0\}
\]
is finite or countable. We index $\Ec$ by letting $I$ be a set and
$I\ni i\mapsto t(i)\in[0,\infty)$ be an injective map such that
\[
\Ec=\{0\}\cup\{E_{|T|}(\{\mu_{t(i)}(T)\})\mid i\in I\},
\]
$E_{|T|}(\{\mu_{t(i)}(T)\})\ne0$ and
\[
t(i)=\inf\{s\mid\mu_s(T)=\mu_{t(i)}(T)\}.
\]
Let $a_i=\tau(E_{|T|}(\{\mu_{t(i)}(T)\}))$.

Fix $i\in I$. Applying the spectral theorem to the normal operator
$TE_{|T|}(\{\mu_{t(i)}(T)\})$ and putting an atomless resolution
of the identity under any of its atoms, we find a family
$(Q_r)_{0\le r<a_i}$ of projections in $\Mcal$ such that
\renewcommand{\labelenumi}{(\arabic{enumi})}
\begin{enumerate}
\item $r_1\le r_2$ implies $Q_{r_1}\le Q_{r_2}$,
\item $\tau(Q_r)=r$,
\item $Q_rT=TQ_r$,
\item $Q_r\le E_{|T|}(\{\mu_{t(i)}(T)\})$.
\end{enumerate}
Let
\[
P_{t(i)+r}=E_{|T|}(\mu_{t(i)}(T),\infty)+Q_r
\]
for all $r\in[0,a_i)$. If $a_i\ne\infty$ then set
$P_{t(i)+a_i}=E_{|T|}[\mu_{t(i)}(T),\infty)$. Now it is easily
seen that the family $(P_t)_{t\ge0}$ satisfies (i)--(v).

Suppose $\lim_{t\to\infty}\mu_t(T)=0$ and let
\[
S=\int_{(0,\infty)}\mu_t(T)\dif F(t).
\]
Clearly $S\ge0$. In order to show $S=|T|$, it will suffice to show
$E_S(x,\infty)=E_{|T|}(x,\infty)$ for all $x>0$. We have
\[
E_S(x,\infty)=F(\{t>0\mid\mu_t(T)>x\}).
\]
But $\{t>0\mid\mu_t(T)>x\}=(0,y)$ where
\[
y=\sup\{t>0\mid\mu_t(T)>x\}=\inf\{t>0\mid\mu_t(T)\le x\}.
\]
From Lemma~\ref{lem:infsup}, $y=\tau(E_{|T|}(x,\infty))$ and,
furthermore, $\mu_y(T)\le x$. By construction,
\[
F((0,y))=P_y=E_{|T|}(x,\infty).
\]
\end{proof}

\section{Classification of modules of a type II factor}
\label{sec:class}

Let $D^+(0,\infty)$, respectively $D^+(0,1)$, denote the cone of
all decreasing (i.e.\ nonincreasing) functions $f$ from the
interval $(0,\infty)$, respectively $(0,1)$, into $[0,\infty)$
that are continuous from the right.

\begin{defi}\rm
Let $\DEu$ be either $D^+(0,\infty)$ or $D^+(0,1)$. A subset
$\Lambda$ of $\DEu$ is called a {\em hereditary subcone} of $\DEu$
if
\renewcommand{\labelenumi}{(\roman{enumi})}
\begin{enumerate}
\item $f,g\in\Lambda$ implies $f+g\in\Lambda$,
\item $f\in\Lambda$, $g\in D^+(0,\infty)$ and $g\le f$ imply $g\in\Lambda$.
\end{enumerate}
The subset $\Lambda\subseteq\DEu$ is called a {\em characteristic
set in} $\DEu$ if it is a hereditary subcone and if
\renewcommand{\labelenumi}{(\roman{enumi})}
\begin{enumerate}
\setcounter{enumi}{2}
\item $f\in\Lambda$ implies $D_2f\in\Lambda$,
\end{enumerate}
where $D_2f(t)=f(t/2)$.
\end{defi}

Let $\Mcal$ be either a type II$_\infty$ factor with a fixed
semifinite normal trace $\tau$ or a type II$_1$ factor with
tracial state $\tau$. Let $\DEu$ be $D^+(0,\infty)$ if $\Mcal$ is
type II$_\infty$ and $D^+(0,1)$ if $\Mcal$ is type II$_1$. We will
recall from~\cite{GuidoIsola} the classification of submodules of
the algebra $\Mcalb$ of $\tau$--measureable operators in terms of
characteristic sets in $\DEu$.

For $T\in\Mcalb$, let $\mu(T)\in\DEu$ be the function which at $t$
takes the value $\mu_t(T)$ of the $t$-th singular number of $T$.
Given a submodule $\Ic\subseteq\Mcalb$, let
\[
\mu(\Ic)=\{\mu(T):T\in\Ic\}\subseteq\DEu.
\]

\begin{prop}[\cite{GuidoIsola}]\label{prop:Class}
Let $\Mcal$ be a factor of type II$_\infty$ or II$_1$. Then the
map $\Ic\mapsto\mu(\Ic)$ is a bijection from the set of all
submodules of $\Mcalb$ onto the set of all characteristic sets in
$\DEu$.
\end{prop}

\begin{remark}\label{rem:obs}\rm
A few well known observations are perhaps in order. If $\Mcal$ is
type II$_\infty$, then $\mu(\Mcal)$ is the set of all bounded
functions in $D^+(0,\infty)$. Thus the smallest nonzero ideal of
$\Mcal$ is the set $\Fc$ of all $\tau$--finite rank
operators in $\Mcal$, where an operator $T$ has $\tau$--finite
rank if $T=ET$ for some projection $E\subseteq\Mcal$ with
$\tau(E)<\infty$; the largest proper ideal of $\Mcal$ is the set
$\Kc$ of all $\tau$--compact operators in $\Mcal$, where
(cf~\cite{StrohWest}) an operator $T$ is $\tau$--compact if
$\lim_{t\to\infty}\mu_t(T)=0$.
\end{remark}

On the other hand, if $\Mcal$ is type II$_1$, then $\mu(\Mcal)$ is
the set of all bounded functions in $D^+(0,1)$, and $\Mcal$ itself
has no proper nonzero ideals.

\section{Sums of commutators}
\label{sec:comm}

\begin{lemma}
\label{lem:nec} Let $\Mcal$ be a von Neumann algebra without
minimal projections and with a normal faithful semifinite trace
$\tau$. Let $T\in\Mcalb$. If $T=\sum_{i=1}^N[A_i,B_i]$ with
$A_i,B_i\in\Mcalb$, then
\begin{equation}
\label{eq:nec}
\big|\tau\big(TE_{|T|}(\mu_s(T),\mu_r(T)]\big)\big|\le
rh(r)+sh(s)
\end{equation}
whenever $0<r<s<\infty$, where
\[
h(t)=(8N+2)\mu_t(T)+(16N+4)\sum_{i=1}^N\mu_t(A_i)\mu_t(B_i).
\]
\end{lemma}
\begin{proof}
We have $E_{|T|}(\mu_s(T),\mu_r(T)]=F_s-F_r$ where
$F_t=E_{|T|}(\mu_t(T),\infty)$. Note
\begin{alignat*}{2}
\|T(1-F_s)\|&\le\mu_s(T),\qquad&\tau(F_s)&\le s, \\
\|T(1-F_r)\|&\le\mu_r(T),\qquad&\tau(F_r)&\le r.
\end{alignat*}
We can find a projection $P\ge F_r$ in $\Mcal$ such that
$\tau(P)\le(4N+1)r$ and
\[
\begin{aligned}
\|A_i(1-P)\|&\le\mu_r(A_i), \\[0.3ex]
\|B_i(1-P)\|&\le\mu_r(B_i),
\end{aligned}
\qquad
\begin{aligned}
\|(1-P)A_i\|&\le\mu_r(A_i), \\[0.3ex]
\|(1-P)B_i\|&\le\mu_r(B_i),
\end{aligned}
\qquad(1\le i\le N).
\]
Then we can find a projection $Q\ge F_s\vee P$ such that
$\tau(Q)\le(4N+1)(r+s)\le(8N+2)s$ and such that for all
$i\in\{1,\ldots,N\}$,
\[
\begin{aligned}
\|A_i(1-Q)\|&\le\mu_s(A_i), \\[0.3ex]
\|B_i(1-Q)\|&\le\mu_s(B_i),
\end{aligned}
\qquad
\begin{aligned}
\|(1-Q)A_i\|&\le\mu_s(A_i), \\[0.3ex]
\|(1-Q)B_i\|&\le\mu_s(B_i),
\end{aligned}
\qquad(1\le i\le N).
\]
Hence
\begin{align*}
|\tau(T(F_s-F_r))|&\le|\tau(T(Q-P))|+|\tau(T(Q-F_s))|+|\tau(T(P-F_r))| \\
&\le|\tau(T(Q-P))|+(8N+2)s\mu_s(T)+(4N+1)r\mu_r(T).
\end{align*}
Since $Q-P$ is a finite projection and $T(Q-P)$ is bounded,
\[
|\tau(T(Q-P))|=|\tau((Q-P)T(Q-P))|\le\sum_{i=1}^N|\tau((Q-P)[A_i,B_i](Q-P))|.
\]
Since also $A_i(Q-P)$, $(Q-P)A_i$, $B_i(Q-P)$ and $(Q-P)B_i$ are
bounded, we have
\begin{align*}
\tau((Q-P)[A_i,B_i](Q-P))&=\tau((Q-P)A_i(1-Q+P)B_i(Q-P)) \\
&-\tau((Q-P)B_i(1-Q+P)A_i(Q-P)).
\end{align*}
But
\begin{align*}
|\tau(&(Q-P)A_i(1-Q+P)B_i(Q-P))| \\
&\le|\tau((Q-P)A_i(1-Q)B_i)|+|\tau(B_i(Q-P)A_iP)| \\
&\le\tau(Q-P)\|A_i(1-Q)\|\,\|(1-Q)B_i\|+\tau(P)\|B_i(Q-P)\|\,\|(Q-P)A_i\| \\
&\le(8N+2)s\mu_s(A_i)\mu_s(B_i)+(4N+1)r\mu_r(B_i)\mu_r(A_i),
\end{align*}
and also with $A_i$ and $B_i$ interchanged. Adding these several
upper bounds gives~\eqref{eq:nec}.
\end{proof}

\begin{lemma}
\label{lem:suf} Let $\Mcal$ be a II$_\infty$ factor with a
specified normal faithful semifinite trace $\tau$. Let $h\in
D^+(0,\infty)$ and suppose $T\in\Mcalb$ is a normal operator
satisfying
\begin{equation}
\label{eq:lim0} \lim_{t\to\infty}\mu_t(T)=0
\end{equation}
and
\begin{equation}
\label{eq:suf}
\big|\tau\big(TE_{|T|}(\mu_s(T),\mu_r(T)]\big)\big|\le
rh(r)+sh(s)\qquad(0<r<s<\infty).
\end{equation}
Let $\Lambda$ be the characteristic subset in $D^+(0,\infty)$
generated by $h$ and $\mu(T)$. Then there are
$X_1,\ldots,X_{14}\in\Mcalb$ and $Y_1,\ldots,Y_{14}\in\Mcal$ such
that
\begin{equation}
\label{eq:T14com} T=\sum_{i=1}^{14}[X_i,Y_i]
\end{equation}
and $\mu(X_i)\in\Lambda$ for all $i$.
\end{lemma}
\begin{proof}
Let $(P_t)_{t\ge0}$ be a family of projections obtained from
Lemma~\ref{lem:mashatoms}. Assumption~\eqref{eq:lim0} implies
$P_\infty T=TP_\infty=T$, where $P_\infty=\bigvee_{t\ge0}P_t$. Let
$P[s,t]=P_t-P_s$ when $s<t$, let
\[
\alpha_n=2^{-n}\tau(TP[2^n,2^{n+1}])\qquad(n\in\Ints)
\]
and let
\[
A=\sum_{n\in\Ints}\alpha_nP[2^n,2^{n+1}].
\]
Note that $T-A=\sum_{n\in\Ints}S_n$ where
$S_n=(T-A)P[2^n,2^{n+1}]$ is an element of the II$_1$--factor
$P[2^n,2^{n+1}]\Mcal P[2^n,2^{n+1}]$ having zero trace.
By~\cite[Thm.\ 2.3]{FackdelaHarpe}, there are
$X_1^{(n)},\ldots,X_{10}^{(n)},Y_1^{(n)},\ldots,Y_{10}^{(n)}\in
P[2^n,2^{n+1}]\Mcal P[2^n,2^{n+1}]$ such that
\[
S_n=\sum_{i=1}^{10}[X_i^{(n)},Y_i^{(n)}]
\]
and for all $i$, $\|X_i^{(n)}\|\le12\|S_n\|$ and
$\|Y_i^{(n)}\|\le2$. We therefore have
\[
T-A=\sum_{i=1}^{10}[X_i,Y_i],
\]
where
\[
X_i=\sum_{n\in\Ints}X_i^{(n)},\qquad
Y_i=\sum_{n\in\Ints}Y_i^{(n)}.
\]
Clearly $Y_i\in\Mcal$. Since
$\|X_i^{(n)}\|\le12\|T_n\|\le12\mu_{2^n}(T)$, it follows that
$\mu_{2^{n+1}}(X_i)\le12\mu_{2^n}(T)$, and therefore that
$\mu(X_i)\in\Lambda$.

It remains to show that $A$ is a sum of four commutators. For
$t>0$ let $F_t=E_{|T|}(\mu_t(T),\infty)$. For $k,\ell\in\Ints$,
$k<\ell$, using the hypothesis~\eqref{eq:suf} we get
\begin{align*}
\bigg|\sum_{j=k}^{\ell-1}2^j\alpha_j\bigg|
&=|\tau(T(P_{2^\ell}-P_{2^k}))| \\
&\le|\tau(T(F_{2^\ell}-F_{2^k}))|+|\tau(T(P_{2^k}-F_{2^k}))|+|\tau(T(P_{2^\ell}-F_{2^\ell}))| \\[1ex]
&\le2^kh(2^k)+2^\ell
h(2^\ell)+2^k\mu_{2^k}(T)+2^\ell\mu_{2^\ell}(T).
\end{align*}
Letting $\phi(t)=h(t)+\mu_t(T)$, we have $\phi\in\Lambda$ and
\begin{equation}
\label{eq:phi}
\bigg|\sum_{j=k}^{\ell-1}2^j\alpha_j\bigg|\le2^k\phi(2^k)+2^\ell\phi(2^\ell).
\end{equation}

We will now write $\RealPart A$ as a sum of two commutators. Note
that inequality~\eqref{eq:phi} continues to hold when each
$\alpha_j$ is replaced by $\RealPart\alpha_j$. We will find real
numbers $\beta_n$ satisfying
\begin{equation}
\label{eq:betan}
\RealPart\alpha_n=\beta_{n-1}-2\beta_n,\qquad|\beta_n|\le\phi(2^n),\qquad(n\in\Ints).
\end{equation}
Treating $\beta_0$ as the independent variable, solving the
equality in~\eqref{eq:betan} recursively yields
\[
\begin{aligned}
\beta_{-m}&=2^m\beta_0+2^{m-1}\sum_{j=-m+1}^02^j\RealPart\alpha_j \\
\beta_m&=2^{-m}\beta_0-2^{-m-1}\sum_{j=1}^m2^j\RealPart\alpha_j
\end{aligned}
\qquad(m\ge1).
\]
The condition $|\beta_n|\le\phi(2^n)$ for all $n\in\Ints$ is thus
equivalent to the inequalitities
\begin{align*}
-2^{-m}\phi(2^{-m})-\frac12\sum_{j=-m+1}^02^j\RealPart\alpha_j&\quad\le\quad\beta_0\quad\le\quad2^{-m}\phi(2^{-m})-\frac12\sum_{j=-m+1}^02^j\RealPart\alpha_j \\
-2^m\phi(2^m)+\frac12\sum_{j=1}^m2^j\RealPart\alpha_j&\quad\le\quad\beta_0\quad\le\quad2^m\phi(2^m)+\frac12\sum_{j=1}^m2^j\RealPart\alpha_j
\end{align*}
for all $m\in\Nats$. The existence of a real number $\beta_0$
satisfying all of these relations is equivalent to the following
four inequalities holding for all integers $k,\ell\ge1$:
\begin{align}
-2^{-k}\phi(2^{-k})-\frac12\sum_{j=-k+1}^02^j\RealPart\alpha_j\quad&\le\quad2^{-\ell}\phi(2^{-\ell})-\frac12\sum_{j=-\ell+1}^02^j\RealPart\alpha_j \label{eq:--} \\
-2^{-k}\phi(2^{-k})-\frac12\sum_{j=-k+1}^02^j\RealPart\alpha_j\quad&\le\quad2^\ell\phi(2^\ell)+\frac12\sum_{j=1}^\ell2^j\RealPart\alpha_j \label{eq:-+} \\
-2^k\phi(2^k)+\frac12\sum_{j=1}^k2^j\RealPart\alpha_j\quad&\le\quad2^{-\ell}\phi(2^{-\ell})-\frac12\sum_{j=-\ell+1}^02^j\RealPart\alpha_j \label{eq:+-} \\
-2^k\phi(2^k)+\frac12\sum_{j=1}^k2^j\RealPart\alpha_j\quad&\le\quad2^\ell\phi(2^\ell)+\frac12\sum_{j=1}^\ell2^j\RealPart\alpha_j.
\label{eq:++}
\end{align}
But these inequalities are easily verified. For
example,~\eqref{eq:--} is equivalent to
\begin{align}
\frac12\sum_{j=-\ell+1}^{-k}2^j\RealPart\alpha_j\quad\le\quad2^{-\ell}\phi(2^{-\ell})+2^{-k}\phi(2^{-k})&\qquad\text{if }k<\ell \label{eq:--1} \\
-2^{-k}\phi(2^{-k})-2^{-\ell}\phi(2^{-\ell})\quad\le\quad\frac12\sum_{j=-k+1}^{-\ell}2^j\RealPart\alpha_j&\qquad\text{if
}k>\ell,
\end{align}
while~\eqref{eq:-+} is equivalent to
\begin{equation}
\label{eq:-+1}
-2^{-k}\phi(2^{-k})-2^\ell\phi(2^\ell)\le\frac12\sum_{j=-k+1}^\ell2^j\RealPart\alpha_j;
\end{equation}
keeping in mind that $\phi$ is nonnegative and nonincreasing,
inequalities~\eqref{eq:--1}--\eqref{eq:-+1} follow directly
from~\eqref{eq:phi}. Inequalities~\eqref{eq:+-} and~\eqref{eq:++}
are verified for all $k$ and $\ell$ similarly. We have suceeded in
proving the existence of $\beta_n$ satisfying~\eqref{eq:betan}.

Now let $V_n,W_n\in\Mcal$, ($n\in\Ints$), be such that
\begin{alignat*}{2}
V_n^*V_n&\quad=\quad P[2^{n-1},2^n],&\qquad V_nV_n^*&\quad=\quad P[2^n,2^n+2^{n-1}], \\
W_n^*W_n&\quad=\quad P[2^{n-1},2^n],&\qquad W_nW_n^*&\quad=\quad
P[2^n+2^{n-1},2^{n+1}]
\end{alignat*}
and let
\begin{alignat*}{2}
X_{11}&=\sum_{n\in\Ints}\beta_{n-1}V_n&\qquad Y_{11}=\sum_{n\in\Ints}V_n \\
X_{12}&=\sum_{n\in\Ints}\beta_{n-1}W_n&\qquad
Y_{12}=\sum_{n\in\Ints}W_n.
\end{alignat*}
Then $X_i\in\Mcalb$, $\mu(X_i)\in\Lambda$ and $Y_i\in\Mcal$
($i=11,12$), and
\[
[X_{11},Y_{11}]+[X_{12},Y_{12}]=\RealPart A.
\]

We may do the same for $\ImagPart A$.
\end{proof}

We now prove an analogous result in a II$_1$--factor.
\begin{lemma}
\label{lem:II1suf} Let $\Mcal$ be a II$_1$--factor with tracial
state $\tau$ and let $T\in\Mcalb$ be a normal operator. Suppose
there is $h\in D^+(0,1)$ such that
\begin{equation}
\label{eq:II1suf}
\big|\tau\big(TE_{|T|}[0,\mu_r(T)]\big)\big|\le
rh(r),\qquad(0<r\le1).
\end{equation}
Let $\Lambda$ be the characteristic set in $D^+(0,1)$ generated by
$h$ and $\mu(T)$. Then there are $X_1,\ldots,X_{12}\in\Mcalb$ and
$Y_1,\ldots,Y_{12}\in\Mcal$ such that
\begin{equation}
\label{eq:II1T14com} T=\sum_{i=1}^{12}[X_i,Y_i]
\end{equation}
and $\mu(X_i)\in\Lambda$ for all $i$.
\end{lemma}
\begin{proof}
Lemma~\ref{lem:mashatoms} (formally applied in $\Mcal\otimes
B(\HEu)$, if we like) gives a family of projections $(P_t)_{0\le
t\le1}$ satisfying (i)--(v) of that proposition. Let
$P[s,t]=P_t-P_s$ ($s<t$), let
\[
\alpha_n=\tau(TP[2^n,2^{n+1}])\qquad(n\in\Ints,\,n<0)
\]
and let
\[
A=\sum_{n=-\infty}^{-1}\alpha_nP[2^n,2^{n+1}].
\]
Applying the result of Fack and de la Harpe~\cite{FackdelaHarpe}
as in the proof of Lemma~\ref{lem:suf}, we can show
\[
T-A=\sum_{i=1}^{10}[X_i,Y_i]
\]
with $Y_i\in\Mcal$, $X_i\in\Mcalb$ and  $\mu(X_i)\in\Lambda$.
Letting $F_t=E_{|T|}(\mu_t(T),\infty)$ and using the
hypothesis~\eqref{eq:II1suf}, for $n\in\Ints$, $n\le-1$ we have
\begin{align*}
\bigg|\sum_{j=n}^{-1}2^j\alpha_j\bigg|&=|\tau(T(1-P_{2^n}))| \\
&\le|\tau(T(1-F_{2^n}))|+|\tau(T(P_{2^n}-F_{2^n}))| \\
&\le2^nh(2^n)+2^n\mu_{2^n}(T).
\end{align*}
Let $\beta_{-1}=0$ and
\[
\beta_n=2^{-n-1}\sum_{j=n+1}^{-1}2^j\alpha_j\qquad(n\in\Ints,\,n\le-2).
\]
Then we have
\[
|\beta_n|\le\frac12(h(2^{n+1})+\mu_{2^{n+1}}(T))
\]
and
\[
\beta_{n-1}-2\beta_n=\alpha_n,\qquad(n\le-1).
\]
Let $V_n,W_n\in\Mcal$ ($n\le-1$) be as in the proof of
Lemma~\ref{lem:suf} and let
\begin{alignat*}{2}
X_{11}&=\sum_{n=-\infty}^{-1}\beta_{n-1}V_n&\qquad Y_{11}=\sum_{n=-\infty}^{-1}V_n \\
X_{12}&=\sum_{n=-\infty}^{-1}\beta_{n-1}W_n&\qquad
Y_{12}=\sum_{n=-\infty}^{-1}W_n.
\end{alignat*}
Then $X_i\in\Mcalb$, $\mu(X_i)\in\Lambda$ and $Y_i\in\Mcal$
($i=11,12$), and
\[
[X_{11},Y_{11}]+[X_{12},Y_{12}]=A.
\]
\end{proof}

\begin{prop}
\label{prop:MM} If $\Mcal$ is a II$_\infty$--factor, then
$[\Mcal,\Mcal]=\Mcal$.
\end{prop}
\begin{proof}
Let $Q_0,\,Q_1,\,\ldots$ be projections in $\Mcal$, each
equivalent to $1$, and such that $\sum_{j=0}^\infty Q_j=1$. Let
$V\in\Mcal$ be such that
\[
V^*V=1,\qquad VV^*=Q_0,
\]
and let $W\in\Mcal$ be such that
\[
W^*W=1,\qquad WQ_jW^*=Q_{j+1}\qquad(j\ge0).
\]
Given $T\in\Mcal$, let
\[
S=\sum_{k=0}^\infty W^kVTV^*(W^*)^k.
\]
Then $S\in\Mcal$ and $[WS,W^*]+[TV^*,V]=T$.
\end{proof}

\begin{lemma}
\label{lem:IJM} Let $\Mcal$ be a II$_1$-- or a
II$_\infty$--factor, and let $\Ic\subseteq\Mcalb$ and
$\Jc\subseteq\Mcalb$ be submodules. Then
$[\Ic\Jc,\Mcal]\subseteq[\Ic,\Jc]$.
\end{lemma}
\begin{proof}
If $X\in\Ic\Jc$ then $X=AB$ for $A\in\Ic$ and $B\in\Jc$. This can
be seen by writing $X=V|X|$ for a partial isometry $V$.

If also $Y\in\Mcal$, then we have
\[
[X,Y]=ABY-YAB=[A,BY]+[B,YA]\in[\Ic,\Jc].
\]
\end{proof}

\begin{thm}
\label{thm:II1} Let $\Mcal$ be a type II$_1$ factor and let
$\Ic\subseteq\Mcalb$ and $\Jc\subseteq\Mcalb$ be submodules. Let
$T\in\Ic\Jc$ be normal. Then $T\in[\Ic,\Jc]$ if and only if there
is $h\in\mu(\Ic\Jc)$ such that
\begin{equation}
\label{eq:necsufII1}
\big|\tau\big(TE_{|T|}[0,\mu_r(T)]\big)\big|\le
rh(r)\qquad(0<r<1).
\end{equation}
\end{thm}
\begin{proof}
We may embed $\Mcal$ in the II$_\infty$ factor $\Mcal\otimes
B(\HEu)$ in a trace--preserving manner. If $T\in[\Ic,\Jc]$ then
letting $h$ be as in Lemma~\ref{lem:nec}, we have
$h\in\mu(\Ic\Jc)$ with $h(1)=0$. So taking $s=1$ in
equation~\eqref{eq:nec}, we get that $h$
satisfies~\eqref{eq:necsufII1}.

Now suppose $h\in\mu(\Ic\Jc)$ is such that~\eqref{eq:necsufII1}
holds. By Lemma~\ref{lem:II1suf}, $T\in[\Ic\Jc,\Mcal]$. Now
Lemma~\ref{lem:IJM} gives $T\in[\Ic,\Jc]$.
\end{proof}

\begin{thm}\label{thm:IIinf}
Let $\Mcal$ be a type II$_\infty$ factor and let
$\Ic\subseteq\Mcalb$ and $\Jc\subseteq\Mcalb$ be submodules. Let
$T\in\Ic\Jc$ be normal. Then $T\in[\Ic,\Jc]$ if and only if there
is $h\in\mu(\Ic\Jc)$ such that
\begin{equation}
\label{eq:necsuf}
\big|\tau\big(TE_{|T|}(\mu_s(T),\mu_r(T)]\big)\big|\le
rh(r)+sh(s)\qquad(0<r<s<\infty).
\end{equation}
\end{thm}
\begin{proof}
If $T\in[\Ic,\Jc]$ then by Lemma~\ref{lem:nec} there is
$h\in\mu(\Ic\Jc)$ satisfying~\eqref{eq:necsuf}.

Now suppose $h\in\mu(\Ic\Jc)$ is such that~\eqref{eq:necsuf}
holds, and let us show $T\in[\Ic\Jc,\Mcal]$. If
$\lim_{t\to\infty}\mu_t(T)=0$, then by Lemma~\ref{lem:suf} we have
$T\in[\Ic\Jc,\Mcal]$. Suppose $d:=\lim_{t\to\infty}\mu_t(T)>0$. If
$T$ is bounded, then by Proposition~\ref{prop:MM},
$T\in[\Mcal,\Mcal]\subseteq[\Ic\Jc,\Mcal]$. Suppose $T$ is
unbounded, let $a>0$ be such that $\mu_a(T)>d$ and let
$Q=E_{|T|}(\mu_a(T),\infty)$. Then $0<\tau(Q)\le a$, $QT=TQ$ and
$\|(1-Q)T\|\le\mu_a(T)$. By Proposition~\ref{prop:MM},
$(1-Q)T\in[\Ic\Jc,\Mcal]$. We have
\[
\mu_t(QT)=\begin{cases}
\mu_t(T)&\text{if }t<\tau(Q) \\
0&\text{if }t\ge\tau(Q).
\end{cases}
\]
Let $0<r<s<\infty$. Then
\[
(QT)E_{|QT|}(\mu_s(QT),\mu_r(QT)]=\begin{cases}
TE_{|T|}(\mu_s(T),\mu_r(T)]&\text{if }s<\tau(Q) \\
TE_{|T|}(\mu_a(T),\mu_r(T)]&\text{if }r<\tau(Q)\le s \\
0&\text{if }\tau(Q)\le r.
\end{cases}
\]
If $r<\tau(Q)\le s$ then we have
\begin{align*}
\big|\tau\big(TE_{|T|}(\mu_a(T),\mu_r(T)]\big)\big|
\le rh(r)+ah(a)&\le rh(r)+ah(\tau(Q)) \\
&\le rh(r)+s\frac{ah(\tau(Q))}{\tau(Q)}.
\end{align*}
Let $\htil(t)=\max\big(h(t),\frac{ah(\tau(Q))}{\tau(Q)}\big)$.
Then $\htil\in\mu(\Ic\Jc)$ and we have
\[
\big|\tau\big((QT)E_{|QT|}(\mu_s(QT),\mu_r(QT)]\big)\big|\le
r\htil(r)+s\htil(s)\qquad(0<r<s<\infty).
\]
Now Lemma~\ref{lem:suf} implies $QT\in[\Ic\Jc,\Mcal]$.

We have shown $T\in[\Ic\Jc,\Mcal]$.
From Lemma~\ref{lem:IJM}, it
follows that $T$ belongs to $[\Ic,\Jc]$.
\end{proof}

\begin{cor}
Let $\Mcal$ be a type II$_\infty$ factor or a type II$_1$ factor,
and let $\Ic\subseteq\Mcalb$ and $\Jc\subseteq\Mcalb$ be
submodules. Then
\[
[\Ic,\Jc]=[\Ic\Jc,\Mcal].
\]
\end{cor}

\section{Separated asymptotic behaviour}
\label{sec:asymp}

Throughout this section, $\Mcal$ will be a type II$_\infty$ factor with semifinite trace $\tau$
and $\Ic\subseteq\Mcalb$ will be a nonzero submodule.
Theorem~\ref{thm:IIinf} gives a necessary and sufficient condition
for a normal operator $T$ to belong to the commutator space
$[\Ic,\Mcal]$, but this condition considers simultaneous
asymptotics at $0$ and $\infty$. In this section, we give an
equivalent characterization which separates the behaviour at $0$
and $\infty$.

We have
\begin{equation}\label{eq:Ifsb}
\Ic=\Ic\fs+\Ic\subb
\end{equation}
where
\begin{align*}
\Ic\fs&=\{T\in\Ic\mid\mu_s(T)=0\text{ for some }s>0\} \\
\Ic\subb&=\{T\in\Ic\mid\mu(T)\text{ bounded }\}.
\end{align*}
Thus $\Ic\fs$ is the set of $T\in\Ic$ that are supported on finite
projections and $\Ic\subb=\Ic\cap\Mcal$. From~\eqref{eq:Ifsb}, we
have
\begin{equation}\label{eq:Icomfsb}
[\Ic,\Mcal]=[\Ic\fs,\Mcal]+[\Ic\subb,\Mcal]\;.
\end{equation}

Given a normal element $T\in\Ic$, using a spectral projection of
$|T|$ we can easily write $T=T\fs+T\subb$ for some normal elements
$T\fs\in\Ic\fs$ and $T\subb\in\Ic\subb$. It is our purpose to use
Theorem~\ref{thm:IIinf} to give necessary and sufficient
conditions for $T\in[\Ic,\Mcal]$ in terms of $T\fs$ and $T\subb$.

\begin{lemma}\label{lem:Ifsb}
Let $\Ic\subseteq\Mcalb$ be a submodule.
\renewcommand{\labelenumi}{(\roman{enumi})}
\begin{enumerate}
\item Let $T\in\Ic\fs$ be normal.
Then $T\in[\Ic\fs,\Mcal]$ if and only if there is
$h\in\mu(\Ic\fs)$ such that
\begin{equation}\label{eq:Trhr}
|\tau(TE_{|T|}[0,\mu_r(T)])|\le rh(r),\quad(0<r<\infty)\;.
\end{equation}
\item
Let $T\in\Ic\subb$ be normal. Then $T\in[\Ic\subb,\Mcal]$ if and
only if there is $h\in\mu(\Ic\subb)$ such that
\begin{equation}\label{eq:Tshs}
|\tau(TE_{|T|}(\mu_s(T),\infty))|\le sh(s),\quad(0<s<\infty)\;.
\end{equation}
\end{enumerate}
\end{lemma}
\begin{proof}
Let us prove~(i). If $T\in[\Ic\fs,\Mcal]$, then invoking
Theorem~\ref{thm:IIinf} and letting $s\to\infty$, since $\mu_s(T)$
and $h(s)$ are eventually zero we obtain
\[
|\tau(TE_{|T|}(0,\mu_r(T)])|\le rh(r)\;,
\]
which clearly implies~\eqref{eq:Trhr}. On the other hand,
if~\eqref{eq:Trhr} holds, then for $0<r<s<\infty$ we have
\begin{align*}
|\tau(TE_{|T|}(\mu_s(T),\mu_r(T)])|
&\le|\tau(TE_{|T|}[0,\mu_r(T)])|+|\tau(TE_{|T|}[0,\mu_s(T)])| \\
&\le rh(r)+sh(s)\;,
\end{align*}
so $T\in[\Ic\fs,\Mcal]$ by Theorem~\ref{thm:IIinf}.

For~(ii), if $T\in[\Ic\subb,\Mcal]$, then invoking
Theorem~\ref{thm:IIinf} and letting $r\to0$, we get
\[
|\tau(TE_{|T|}(\mu_s(T),\infty))|=|\tau(TE_{|T|}(\mu_s(T),\|T\|])|\le
sh(s)\;,
\]
since $h(r)$ stays bounded. Thus~\eqref{eq:Tshs} holds. The
argument that~\eqref{eq:Tshs} implies $T\in[\Ic\subb,\Mcal]$ is
similar to the analogous one in case~(i).
\end{proof}

Recall (Remark~\ref{rem:obs})
$\Fc$ denotes the submodule (in fact, the ideal of $\Mcal$) consisting of $\tau$--finite rank bounded operators:
$\Fc=\Mcal\fs$.
\begin{cor}\label{cor:MfsM}
$[\Fc,\Mcal]=\Fc\cap\ker\tau$.
\end{cor}
\begin{proof}
It will suffice to show that if $T=T^*\in\Fc$, then
$T\in[\Fc,\Mcal]$ if and only if $\tau(T)=0$. Suppose
$T\in[\Fc,\Mcal]$ and let $h\in\mu(\Fc)$ be such
that~\eqref{eq:Tshs} holds. Then $h(s)=0$ and $\mu_s(T)=0$ for
some $s>0$ and therefore
$\tau(T)=\tau(TE_{|T|}(\mu_s(T),\infty))=0$.

Suppose $\tau(T)=0$. Then $\mu_{s'}(T)=0$ for some $s'>0$. Let
\[
h(s)=\begin{cases} \|T\|&s<s' \\ 0&s\ge s'\;.\end{cases}
\]
Then $h\in\mu(\Fc)$. Using~\eqref{eq:tauE}, we see
that~\eqref{eq:Tshs} holds when $0<s<s'$, and it holds when $s\ge
s'$ because $\tau(T)=0$.
\end{proof}

See Definition~\ref{defi:oplus} for an explanation of the notation $\oplus$
used below.

\begin{thm}\label{thm:oplusX}
Let $\Ic\subseteq\Mcalb$ be a nonzero submodule and let
$T=T\fs+T\subb\in\Ic$, where $T\fs\in\Ic\fs$ and
$T\subb\in\Ic\subb$. Then the following are equivalent:
\renewcommand{\labelenumi}{(\alph{enumi})}
\begin{enumerate}
\item $T\in[\Ic,\Mcal]$.
\item There is $X\in\Fc$ such that
\begin{equation}\label{eq:TX}
\begin{aligned}
T\subb\oplus X&\in[\Ic\subb,\Mcal] \\
T\fs\oplus(-X)&\in[\Ic\fs,\Mcal]\;.
\end{aligned}
\end{equation}
\item There is $a\in\Cpx$ such that whenever $X,Y\in\Fc$, $\tau(X)\ne0$ and $\tau(Y)\ne0$,
\begin{align}
T\subb\oplus\tfrac a{\tau(X)}X&\in[\Ic\subb,\Mcal] \label{eq:TbaX} \\
T\fs\oplus\tfrac{-a}{\tau(Y)}Y&\in[\Ic\fs,\Mcal]\;.
\label{eq:TfsaX}
\end{align}
\end{enumerate}
\end{thm}
\begin{proof}
We first prove (a)$\implies$(c). Suppose $T\in[\Ic,\Mcal]$.
From~\eqref{eq:Icomfsb}, we have $T=\Tt\fs+\Tt\subb$ for some
$\Tt\fs\in[\Ic\fs,\Mcal]$ and $\Tt\subb\in[\Ic\subb,\Mcal]$. Then
using Corollary~\ref{cor:MfsM},
\[
\Tt\subb-T\subb=T\fs-\Tt\fs\in\Ic\subb\cap\Ic\fs=\Fc\;.
\]
Let $a=\tau(\Tt\subb-T\subb)$ and let $X\in\Fc$ with
$\tau(X)\ne0$. Then
\[
\Tt\subb-T\subb-\tfrac
a{\tau(X)}X\in\Fc\cap\ker\tau=[\Fc,\Mcal]\subseteq[\Ic\subb,\Mcal]\;.
\]
Thus
\begin{align*}
T\subb\oplus\tfrac a{\tau(X)}X&\in(T\subb\oplus\tfrac a{\tau(X)}X\oplus0)+[\Ic\subb,\Mcal] \\
&=(T\subb\oplus\tfrac a{\tau(X)}X\oplus(\Tt\subb-T\subb-\tfrac a{\tau(X)}X))+[\Ic\subb,\Mcal] \\
&=(T\subb\oplus(\Tt\subb-T\subb))+[\Ic\subb,\Mcal] \\
&=\Tt\subb+[\Ic\subb,\Mcal]=[\Ic\subb,\Mcal]
\end{align*}
and~\eqref{eq:TbaX} holds. Similarly, we have
\begin{align*}
T\fs\oplus\tfrac{-a}{\tau(Y)}Y&\in(T\fs\oplus\tfrac{-a}{\tau(Y)}Y\oplus0)+[\Ic\fs,\Mcal] \\
&=(T\fs\oplus\tfrac{-a}{\tau(Y)}Y\oplus(\Tt\fs-T\fs+\tfrac a{\tau(Y)}Y))+[\Ic\fs,\Mcal] \\
&=(T\fs\oplus(\Tt\fs-T\fs))+[\Ic\fs,\Mcal] \\
&=\Tt\fs+[\Ic\fs,\Mcal]=[\Ic\fs,\Mcal]
\end{align*}
and~\eqref{eq:TfsaX} holds.

The implication (c)$\implies$(b) is clear.

For (b)$\implies$(a), assuming~\eqref{eq:TX}, we have
\begin{align*}
T\fs+T\subb&\in T\fs\oplus T\subb+[\Ic,\Mcal] \\
&=T\fs\oplus(-X)\oplus X\oplus T\subb+[\Ic,\Mcal]=[\Ic,\Mcal]\;.
\end{align*}
\end{proof}

\begin{lemma}\label{lem:ataufs}
Let $\Ic\subseteq\Mcalb$ be a nonzero submodule, let $T\in\Ic\fs$
be normal, $T\ne0$ and let $a\in\Cpx$. Let $P\in\Fc$ be a
nonzero projection such that either $T$ is unbounded or
$\frac{|a|}{\tau(P)}<\|T\|$. Then
\begin{equation*}
T\oplus\tfrac a{\tau(P)}P\in[\Ic\fs,\Mcal]
\end{equation*}
if and only if there is $h\in\mu(\Ic\fs)$ such that
\begin{equation}\label{eq:lemafs}
\forall r\in(0,1),\quad |a+\tau(TE_{|T|}[0,\mu_r(T)])|\le
rh(r)\;.
\end{equation}
\end{lemma}
\begin{remark}\rm
As will be apparent from the proof, for any $r'>0$ the existence
of $h\in\mu(\Ic\fs)$ such that~\eqref{eq:lemafs} holds is
equivalent to the existence of $h'\in\mu(\Ic\fs)$ such that
\[
\forall r\in(0,r'),\quad |a+\tau(TE_{|T|}[0,\mu_r(T)])|\le rh'(r)
\]
holds.
\end{remark}
\begin{proof}[Proof Lemma~\ref{lem:ataufs}]
There is $r'>0$ such that $\mu_r(T)>\frac{|a|}{\tau(P)}$ for all
$r\in(0,r')$. Let $T'=T\oplus\frac a{\tau(P)}P$. Then (by
Proposition~\ref{prop:SoplusT}), for $r\in(0,r')$ we have
$\mu_r(T')=\mu_r(T)$,
\begin{align*}
E_{|T'|}[0,\mu_r(T')]&=E_{|T|}[0,\mu_r(T)]\oplus P \\
\tau(T'E_{|T'|}[0,\mu_r(T')])&=a+\tau(TE_{|T|}[0,\mu_r(T)])\;.
\end{align*}
If $T'\in[\Ic\fs,\Mcal]$, then by Lemma~\ref{lem:Ifsb}, there is
$h'\in\mu(\Ic\fs)$ such that
\[
\forall r\in(0,r'),\quad |a+\tau(TE_{|T|}[0,\mu_r(T)])|\le
rh'(r)\;.
\]
Since $d:=\tau(E_{|T|}(0,\infty))<\infty$ and
\[
|\tau(TE_{|T|}[0,\mu_r(T)])|\le\mu_r(T)d
\]
for all $r>0$, we can find $h\in\mu(\Ic\fs)$ such
that~\eqref{eq:lemafs} holds.

Conversely, suppose $h\in\mu(\Ic\fs)$ is such
that~\eqref{eq:lemafs} holds. Assume without loss of generality
$r'\le 1$. Then we have
\[
|\tau(T'E_{|T'|}[0,\mu_r(T')])|\le rh(r)
\]
for all $r\in(0,r')$. Let $r''>r'$ be such that $\mu_{r''}(T')=0$.
Let $d'=E_{|T'|}(0,\infty)$. Then
\[
|\tau(T'E_{|T'|}[0,\mu_r(T')])|\le\begin{cases}
0&\text{if }r\ge r'' \\
\mu_r(T')d'&\text{otherwise.}
\end{cases}
\]
Letting
\[
h'(t)=\begin{cases}
\max(h(t),\frac{\mu_{r'}(T')d'}{r'})&\text{if }0<t<r' \\
\frac{\mu_{r'}(T')d'}{r'}&\text{if }r'\le t<r'' \\
0&\text{if }r''\le t\;,
\end{cases}
\]
we have $h'\in\mu(\Ic\fs)$ and
\[
|\tau(T'E_{|T'|}[0,\mu_r(T')])|\le rh'(r)
\]
for all $r>0$. Thus $T'\in[\Ic\fs,\Mcal]$ by Lemma~\ref{lem:Ifsb}.
\end{proof}

\begin{lemma}\label{lem:ataub}
Let $\Ic\subseteq\Mcalb$ be a nonzero submodule, let
$T\in\Ic\subb$ be normal and let $a\in\Cpx$. If $a\ne0$, let
$P\in\Fc$ be a projection such that $\frac{|a|}{\tau(P)}>\|T\|$.
If $a=0$, let $P\in\Fc$ have nonzero trace. Then
\begin{equation*}
T\oplus\tfrac a{\tau(P)}P\in[\Ic\subb,\Mcal]
\end{equation*}
if and only if there is $h\in\mu(\Ic\subb)$ such that
\begin{equation}\label{eq:lemab}
\forall s\in[1,\infty),\quad
|a+\tau(TE_{|T|}(\mu_s(T),\infty))|\le sh(s)\;.
\end{equation}
\end{lemma}
\begin{remark}\rm
As will be apparent from the proof, for any $s'>0$ the existence
of $h\in\mu(\Ic\subb)$ such that~\eqref{eq:lemab} holds is
equivalent to the existence of $h'\in\mu(\Ic\subb)$ such that
\[
\forall s\in[s',\infty),\quad
|a+\tau(TE_{|T|}(\mu_s(T),\infty))|\le sh'(s)
\]
holds.
\end{remark}
\begin{proof}[Proof Lemma~\ref{lem:ataub}]
Suppose $a\ne0$. Let $T'=T\oplus\frac a{\tau(P)}P$. Then for all
$s>0$, we have, (by Proposition~\ref{prop:SoplusT}),
$\mu_{s+\tau(P)}(T')=\mu_s(T)$,
\begin{align}
E_{|T'|}(\mu_{s+\tau(P)}(T'),\infty)&=E_{|T|}(\mu_{s}(T),\infty)\oplus P\;, \notag \\
\tau(T'E_{|T'|}(\mu_{s+\tau(P)}(T'),\infty))&=a+\tau(TE_{|T|}(\mu_{s}(T),\infty))\;.
\label{eq:tauTprime}
\end{align}
If $T'\in[\Ic\subb,\Mcal]$, then it follows
from~\eqref{eq:tauTprime} and Lemma~\ref{lem:Ifsb} that there is
$h'\in\mu(\Ic\subb)$ such that
\[
\forall s\in(0,\infty),\quad
|a+\tau(TE_{|T|}(\mu_s(T),\infty))|\le(s+\tau(P))h'(s+\tau(P))\;.
\]
Letting $h(s)=(1+\tau(P))h'(s+\tau(p))$, we have
$h\in\mu(\Ic\subb)$ and that~\eqref{eq:lemab} holds.

On the other hand, still taking $a\ne0$, suppose
$h\in\mu(\Ic\subb)$ and~\eqref{eq:lemab} holds.
Using~\eqref{eq:tauTprime}, we have
\[
|\tau(T'E_{|T'|}(\mu_t(T'),\infty))|\le(t-\tau(P))h(t-\tau(P))
\]
for all $t\ge1+\tau(P)$. Using Proposition~\ref{prop:FK}, we have
\[
|\tau(T'E_{|T'|}(\mu_t(T'),\infty))|\le\|T'\|t
\]
for all $t>0$. Therefore, letting
\[
h'(t)=\begin{cases}
\frac1{1+\tau(P)}h(t-\tau(P))&\text{if }t\ge1+\tau(P) \\
\max(\frac{|a|}{\tau(P)},\frac{h(1)}{1+\tau(P)})&\text{if
}0<t<1+\tau(P),
\end{cases}
\]
we get $h'\in\mu(\Ic\subb)$ and
\[
|\tau(T'E_{|T'|}(\mu_t(T'),\infty))|\le th'(t)
\]
for all $t>0$. Thus $T'\in[\Ic\subb,\Mcal]$ by
Lemma~\ref{lem:Ifsb}.

When $a=0$, the existence of $h\in\mu(\Ic\subb)$
satifying~\eqref{eq:lemab} follows from $T\in[\Ic\subb,\Mcal]$
directly from Lemma~\ref{lem:Ifsb}, while proving that the
existence of $h\in\mu(\Ic\subb)$ such that~\eqref{eq:lemab} holds
implies $T\in[\Ic\subb,\Mcal]$ is similar to the case $a\ne0$, but
easier.
\end{proof}

\begin{thm}\label{thm:TfsTbnormal}
Let $\Ic\subseteq\Mcalb$ be a nonzero submodule and let
$T=T\fs+T\subb\in\Ic$, where $T\fs\in\Ic\fs$ and
$T\subb\in\Ic\subb$ are normal. Then $T\in[\Ic,\Mcal]$ if and only
if there are $a\in\Cpx$, $h\fs\in\mu(\Ic\fs)$ and
$h\subb\in\mu(\Ic\subb)$ such that
\begin{align}
\forall r\in(0,1)\quad
|a-\tau(T\fs E_{|T\fs|}[0,\mu_r(T\fs)])|&\le rh\fs(r) \label{eq:afs} \\
\forall s\in[1,\infty)\quad |a+\tau(T\subb
E_{|T\subb|}(\mu_s(T\subb),\infty))|&\le sh\subb(s)\;.
\label{eq:ab}
\end{align}
\end{thm}
\begin{proof}
If $T\fs\ne0$, then the conclusion of the theorem follows from Theorem~\ref{thm:oplusX} and
Lemmas~\ref{lem:ataufs} and~\ref{lem:ataub}.
If $T\fs=0$, then we choose $a=0$ and apply Lemma~\ref{lem:Ifsb}.
\end{proof}

Let $\omega\fs,\omega\subb\in D^+(0,\infty)$ be given by
\begin{align*}
\omega\fs(t)&=\begin{cases}
1/t&\text{if }t<1 \\
0&\text{if }t\ge1,
\end{cases} \\
\omega\subb(t)&=\frac1{1+t}\;.
\end{align*}

\begin{cor}\label{cor:omega}
Let $\Ic\subseteq\Mcalb$ be a nonzero submodule and let
$T=T\fs+T\subb\in\Ic$, where $T\fs\in\Ic\fs$ and
$T\subb\in\Ic\subb$ are normal.
\renewcommand{\labelenumi}{(\Roman{enumi})}
\begin{enumerate}
\item
Suppose $\omega\fs,\omega\subb\in\mu(\Ic)$. Then $T\in[\Ic,\Mcal]$
if and only if $T\fs\in[\Ic\fs,\Mcal]$ and
$T\subb\in[\Ic\subb,\Mcal]$.

\vskip1ex
\item
Suppose $\omega\fs\in\mu(\Ic)$ and $\omega\subb\not\in\mu(\Ic)$.
Then $T\in[\Ic,\Mcal]$ if and only if $T\fs\in[\Ic\fs,\Mcal]$ and
there are $a\in\Cpx$ and $h\subb\in\mu(\Ic\subb)$ such
that~\eqref{eq:ab} holds.

\vskip1ex
\item
Suppose $\omega\fs\not\in\mu(\Ic)$ and $\omega\subb\in\mu(\Ic)$.
Then $T\in[\Ic,\Mcal]$ if and only if $T\subb\in[\Ic\subb,\Mcal]$
and there are $a\in\Cpx$ and $h\fs\in\mu(\Ic\fs)$ such
that~\eqref{eq:afs} holds.
\end{enumerate}
\end{cor}
\begin{proof}
If $\omega\subb\in\mu(\Ic)$, then for any $a\in\Cpx$, the function
\[
t\mapsto\begin{cases}
|a|/t,&0<t<1, \\ 0,&t\ge1
\end{cases}
\]
lies in $\mu(\Ic\fs)$, while if $\omega\fs\in\mu(\Ic)$, then for any $a\in\Cpx$, the function
\[
t\mapsto\begin{cases}
|a|,&t\in(0,1), \\ |a|/t,&t\ge1
\end{cases}
\]
lies in $\mu(\Ic\subb)$.
\end{proof}

This seems like a convenient place to prove the following proposition, which will be needed in Section~\ref{spectral}.
\begin{prop}\label{prop:I0}
Let $\Ic\subseteq\Mcalb$ be a nonzero submodule and suppose $\Mcal\subseteq\Ic$.
Let
\[
\Ic_0=\{T\in\Ic\mid\lim_{t\to\infty}\mu_t(T)=0\}.
\]
Then $[\Ic,\Mcal]\cap\Ic_0=[\Ic_0,\Mcal]$.
\end{prop}
\begin{proof}
Since $\supseteq$ is clear, we need only show $\subseteq$.
Suppose $T\in[\Ic,\Mcal]\cap\Ic_0$ and $T$ is normal.
It will suffice to show $T\in[\Ic_0,\Mcal]$.
Let $T=T\fs+T\subb$ where $T\fs\in(\Ic_0)\fs=\Ic\fs$ and $T\subb\in(\Ic_0)\subb$ are normal.
Note we have $(\Ic_0)\subb=\Mcal_0=\Kc$, (see Remark~\ref{rem:obs}).
Since $T\in[\Ic,\Mcal]$,
by Corollary~\ref{cor:omega}, if $\omega\fs\in\mu(\Ic)$, then $T\fs\in[\Ic\fs,\Mcal]$,
while if $\omega\fs\notin\mu(\Ic)$, then there are $a\in\Cpx$ and $h\fs\in\mu(\Ic\fs)$
such that~\eqref{eq:afs} holds.
Since $\omega\subb\in\mu(\Ic_0)$, by Corollary~\ref{cor:omega} in order to show $T\in[\Ic_0,\Mcal]$
it will suffice to show $T\subb\in[(\Ic_0)\subb,\Mcal]$.
But $[\Kc,\Mcal]=\Kc\ni T\subb$.
\end{proof}

We will finish this section with a few observations relating $[\Ic,\Mcal]$ to $[\Ic_b,\Mcal]$ and $[\Ic\fs,\Mcal]$,
and examples involving ideals of $p$--summable operators.
Writing $\Ic=\Ic\fs+\Ic_b$, we have $[\Ic,\Mcal]=[\Ic\fs,\Mcal]+[\Ic_b,\Mcal]$.
Since $\Ic\fs\cap\Ic_b=\Fc$, and (see Corollary~\ref{cor:MfsM}) $[\Fc,\Mcal]=\Fc\cap\ker\tau$, we have
\begin{align*}
[\Ic\fs,\Mcal]\cap\Ic_b&=[\Ic\fs,\Mcal]\cap\Fc=\begin{cases}
\Fc&\text{if }\omega\fs\in\mu(\Ic) \\
\Fc\cap\ker\tau&\text{if }\omega\fs\notin\mu(\Ic),
\end{cases} \\[1ex]
[\Ic_b,\Mcal]\cap\Ic\fs&=[\Ic_b,\Mcal]\cap\Fc=\begin{cases}
\Fc&\text{if }\omega_b\in\mu(\Ic) \\
\Fc\cap\ker\tau&\text{if }\omega_b\notin\mu(\Ic).
\end{cases}
\end{align*}
So we have the following result.
\begin{prop}\label{prop:IcMasympt}
Let $\Ic$ be a nonzero submodule of $\Mcalb$, for a II$_\infty$ factor $\Mcal$.
Then
\renewcommand{\labelenumi}{(\roman{enumi})}
\begin{enumerate}
\item
$\Fc+[\Ic,\Mcal]=\Ic$ if and only if 
$\Fc+[\Ic\fs,\Mcal]=\Ic\fs$ and
$\Fc+[\Ic_b,\Mcal]=\Ic_b$;
\item
$[\Ic,\Mcal]=\Ic$ if and only if at least one of the following holds:
\renewcommand{\labelenumii}{(\alph{enumii})}
\begin{enumerate}
\item
$[\Ic\fs,\Mcal]=\Ic\fs$ and
$\Fc+[\Ic_b,\Mcal]=\Ic_b$;
\item
$\Fc+[\Ic\fs,\Mcal]=\Ic\fs$ and
$[\Ic_b,\Mcal]=\Ic_b$.
\end{enumerate}
\end{enumerate}
\end{prop}

We now relate the commutator space $[\Ic_b,\Mcal]$ to its discrete analogue.
Let $\Bc\subseteq\Mcal$ be any type I$_\infty$ factor (i.e.\ a copy of $B(\HEu)$)
such that the restriction of $\tau$ to $\Bc$ is semifinite.
Let $\Ic_d=\Ic\cap\Bc$ and let $\Fc_d=\Fc\cap\Bc$; (the ``d'' is for ``discrete'').
Note that $\Ic_d$ is an ideal of $\Bc$ and $\Fc_d$ is the ideal of finite rank operators in $\Bc$.
In the notation used in~\cite{DFWW}, the characteristic set $\mu(\Ic_d)$ of $\Ic_d$,
consisting of the sequences of singular numbers of elements of $\Ic_d$, is naturally identified
with the set of all functions $f\in\mu(\Ic)$ that are constant on the intervals
$[0,1),\,[1,2),\,[2,3),\,\ldots$.
The commutator space $[\Ic_d,\Bc]$ of an ideal of a I$_\infty$ factor has been extensively studied
--- see~\cite{DFWW} and references contained therein, and see~\cite{KW} for some further results.
\begin{lemma}\label{lem:TIb}
Let $T\in\Ic_b$ and assume $\lim_{t\to\infty}\mu_t(T)=0$.
Then there is $A\in\Ic_d$ such that $T-A\in[\Ic_b,\Mcal]$.
\end{lemma}
\begin{proof}
We may without loss of generality assume $T=T^*$ and that $\tau(\Qt)=1$ for a minimal projection $\Qt$ of $\Bc$.
Let $(P_t)_{t\ge0}$ be a family of projections in $\Mcal$ obtained from Lemma~\ref{lem:mashatoms}.
Let $Q_k=P_k-P_{k-1}$, ($k\in\Nats$), $\alpha_k=\tau(TQ_k)$ and $A'=\sum_{k=1}^\infty\alpha_kQ_k$.
Then $TQ_k-\alpha_kQ_k$ is an element of the II$_1$--factor $Q_k\Mcal Q_k$ of trace zero and with
\[
\|(T-\alpha_k)Q_k\|\le\|TQ_k\|+|\alpha_k|\le2\|TQ_k\|\le2\mu_k(T).
\]
Using~\cite[Thm.\ 2.3]{FackdelaHarpe} as in the proof of Lemma~\ref{lem:suf},
one shows $T-A'\in[\Ic_b,\Mcal]$.
Let $\Qt_1,\Qt_2,\ldots\in\Bc$ be pairwise orthogonal projections, each of trace $1$, and let
$U\in\Mcal$ be a partial isometry such that $U^*Q_jU=\Qt_j$.
Let $A=\sum_{k=1}^\infty\alpha_k\Qt_k$.
Then $A=U^*A'U\in\Ic_d$ and $A'-A=[U,U^*A']\in[\Ic_b,\Mcal]$.
Thus $T-A\in[\Ic_b,\Mcal]$.
\end{proof}

\begin{prop}\label{prop:IbId}
\renewcommand{\labelenumi}{(\roman{enumi})}
\begin{enumerate}
\item $\Bc\cap[\Ic_b,\Mcal]=[\Ic_d,\Bc]$.
\item $[\Ic_b,\Mcal]=\Ic_b$ if and only if $[\Ic_d,\Bc]=\Ic_d$.
\item $\Fc+[\Ic_b,\Mcal]=\Ic_b$ if and only if $\Fc_d+[\Ic_d,\Bc]=\Ic_d$.
\end{enumerate}
\end{prop}
\begin{proof}
We may without loss of generality assume $\tau(F)=1$ for a minimal projection $F$ of $\Bc$.
The inclusion $\supseteq$ in~(i) is clear.
To show $\subseteq$, it will suffice to show that $T=T^*\in\Bc\cap[\Ic_b,\Mcal]$ implies $T\in[\Ic_d,\Bc]$.
By Lemma~\ref{lem:Ifsb}, there is $h\in\mu(\Ic_b)$ satisfying~\eqref{eq:Tshs}.
Since $h$ is bounded, replacing $h$ if necessary by a slightly greater function, we may without loss of generality assume
$h$ is constant on all intervals $[0,1),\,[1,2),\,\ldots$.
We may write $T=\sum_{i=1}^\infty\lambda_iF_i$
for a sequence of pairwise orthogonal, minimal projections $F_i$ of $\Bc$ and
for $\lambda_i\in\Reals$ with $|\lambda_1|\ge|\lambda_2|\ge\cdots$.
If $\lim_{n\to\infty}|\lambda_n|>0$, then $\Ic_b=\Mcal$ and $\Ic_d=\Bc$, so~(i) holds.
Hence we may without loss of generality assume $\lim_{n\to\infty}|\lambda_n|=0$.
Suppose $k$ and $n$ are nonnegative integers with $k<n$,
\[
|\lambda_{k+1}|=|\lambda_{k+2}|=\cdots=|\lambda_n|>|\lambda_{n+1}|
\]
and either $k=0$ or $|\lambda_k|>|\lambda_{k+1}|$.
If $s\in[k,n)$, then $\mu_s(T)=|\lambda_{k+1}|$, so by~\eqref{eq:Tshs},
\[
|\lambda_1+\cdots+\lambda_k|=|\tau(TE_{|T|}(\mu_s(T),\infty))|\le sh(s).
\]
Thus, if $\ell\in\{k,\ldots,n-1\}$ and $\ell\ne0$, then
\[
|\lambda_1+\cdots+\lambda_\ell|\le|\lambda_1+\cdots+\lambda_k|+(\ell-k)|\lambda_\ell|\le\ell h(\ell)+\ell|\lambda_\ell|
\]
and
\[
\frac{|\lambda_1+\cdots+\lambda_\ell|}\ell\le h(\ell)+|\lambda_\ell|.
\]
From this, the main result of~\cite{DFWW} implies $T\in[\Ic_d,\Bc]$,
and~(i) is proved.

From~(i), we have
\[
[\Ic_b,\Mcal]=\Ic_b\qquad\implies\qquad[\Ic_d,\Bc]=\Ic_d.
\]
The reverse implication follows from Lemma~\ref{lem:TIb}.
Hence~(ii) is proved.

To prove~(iii), we have $\Fc=\Fc_d+(\Fc\cap\ker\tau)=\Fc_d+[\Fc,\Mcal]$,
so
\[
\Fc+[\Ic_b,\Mcal]=\Fc_d+[\Ic_b,\Mcal].
\]
From~(i) we thus obtain
\[
\Fc+[\Ic_b,\Mcal]=\Ic_b\qquad\implies\qquad\Fc_d+[\Ic_d,\Bc]=\Ic_d.
\]
The reverse implication follows from Lemma~\ref{lem:TIb}.
\end{proof}

We now point out results relating $[\Ic\fs,\Mcal]$ and commutator spaces of submodules of II$_1$--factors.
Let $P\in\Mcal$ be a projection with $\tau(P)=1$ and consider the II$_1$--factor $\Mcal_1=P\Mcal P$.
Then $P\Mcalb P$ is equal to the module $\overline{\Mcal_1}$ of $\tau$--measureable operators affiliated to $\Mcal_1$.
Given a nonzero submodule $\Ic$ of $\Mcalb$, consider the submodule $\Ic_1=P\Ic P$ of $\overline{\Mcal_1}$.
Then the following result follows directly from the characterizations of commutator spaces found in Theorem~\ref{thm:II1}
and Lemma~\ref{lem:Ifsb}.
\begin{prop}\label{prop:IfsII1}
\renewcommand{\labelenumi}{(\roman{enumi})}
\begin{enumerate}
\item $\overline{\Mcal_1}\cap[\Ic\fs,\Mcal]=[\Ic_1,\Mcal_1]$.
\item $[\Ic\fs,\Mcal]=\Ic\fs$ if and only if $[\Ic_1,\Mcal_1]=\Ic_1$.
\item $\Fc+[\Ic\fs,\Mcal]=\Ic\fs$ if and only if $\Mcal_1+[\Ic_1,\Mcal_1]=\Ic_1$.
\end{enumerate}
\end{prop}

For $0<p<\infty$, let $\Lc_p$ denote the submodule of $\Mcalb$ whose characterisitc set $\mu(\Lc_p)$
consists of all the $p$--integrable functions in $D^+(0,\infty)$.
Thus
\[
\Lc_p=\{T\in\Mcalb\mid\tau((T^*T)^{p/2})<\infty\},
\]
where we have extended $\tau$ in the usual way to be a map from positive elements of $\Mcalb$ to $[0,+\infty]$.
Also, let $\Lc_\infty=\Mcal$.

\begin{prop}\label{prop:Lp}
If $0<p<1$, then
\begin{equation}\label{eq:Lp1}
[(\Lc_p)\fs,\Mcal]=(\Lc_p)\fs
\end{equation}
and
\begin{equation}\label{eq:Lp2}
[(\Lc_p)_b,\Mcal]=(\Lc_p)_b\cap\ker\tau,
\end{equation}
so $\Fc+[(\Lc_p)_b,\Mcal]=(\Lc_p)_b$.

With $p=1$, we have
\begin{equation}\label{eq:Lp3}
\Fc+[(\Lc_1)\fs,\Mcal]\ne(\Lc_1)\fs
\end{equation}
and
\begin{equation}\label{eq:Lp4}
\Fc+[(\Lc_1)_b,\Mcal]\ne(\Lc_1)_b.
\end{equation}

If $1<p\le\infty$, then
\begin{equation}\label{eq:Lp5}
[(\Lc_p)\fs,\Mcal]=(\Lc_p)\fs\cap\ker\tau,
\end{equation}
so $\Fc+[(\Lc_p)\fs,\Mcal]=(\Lc_p)\fs$, and
\begin{equation}\label{eq:Lp6}
[(\Lc_p)_b,\Mcal]=(\Lc_p)_b.
\end{equation}
\end{prop}
\begin{proof}
When $p=\infty$, we have $(\Lc_p)\fs=\Fc$ and $(\Lc_p)_b=\Mcal$, and these special cases
of~\eqref{eq:Lp5} and~\eqref{eq:Lp6} have been considered previously.
For $p<\infty$, all of the relations \eqref{eq:Lp1}--\eqref{eq:Lp6} can be readily verified from
properties of $L^p$--functions.

Moreover,~\eqref{eq:Lp2}, \eqref{eq:Lp4} and~\eqref{eq:Lp6} follow from Proposition~\ref{prop:IbId}
and the coresponding discrete analogues, which follow readily from the main result of~\cite{DFWW}
and were originally proved in~\cite{Anderson}, \cite{Weiss} and~\cite{PearcyTopping}, respectively.
On the other hand,~\eqref{eq:Lp1} and~\eqref{eq:Lp5} follow from Proposition~\ref{prop:IfsII1} and~\cite[Prop.\ 2.12]{FigielKalton}.

As an example, let us verify~\eqref{eq:Lp3} directly.
Clearly $[(\Lc_1)\fs,\Mcal]\subseteq\ker\tau$, so it will suffice to find $T=T^*\in(\Lc_1)\fs\cap\ker\tau$
with $T\notin[(\Lc_1)\fs,\Mcal]$.
Using Lemma~\ref{lem:Ifsb}, it will suffice to find $f\in L^1[0,1]$ such that $\int_0^1f=0$ but
the function
\[
s\mapsto\frac1s\int_s^1f(t)dt,\qquad0<s<1
\]
is not integrable.
Such a function is given by
\[
f(t)=\begin{cases}
\frac1{t(\log t)^2}&\text{if }0<t<1/2 \\
\frac{-2}{\log 2}&\text{if }1/2\le t<1.
\end{cases}
\]
\end{proof}

Propositions~\ref{prop:Lp} and~\ref{prop:IcMasympt} now yield the following examples.
\begin{examples}\rm
Let $\Ic=(\Lc_p)\fs+(\Lc_q)_b$, for some $0<p,q\le\infty$.
\renewcommand{\labelenumi}{(\roman{enumi})}
\begin{enumerate}
\item If $p<1$ and $q\ne1$ or if $p\ne1$ and $q>1$, then $[\Ic,\Mcal]=\Ic$.
\item If $p>1$ and $q<1$, then $[\Ic,\Mcal]=\Ic\cap\ker\tau$ and $\Fc+[\Ic,\Mcal]=\Ic$.
\item If $p=1$ or $q=1$, then $\Fc+[\Ic,\Mcal]\ne\Ic$.
\end{enumerate}
\end{examples}

\section{Spectral characterization of $[\mathcal I,\mathcal M].$}\label{spectral}

In this section, $\mathcal M$ will be a II$_{\infty}$--factor with fixed normal, semifinite trace $\tau$.

Let $\mathcal L_{\log}$ be the submodule of all
$T\in\overline{\mathcal M}$ such that
$$ \int_0^{\infty}\log (1+\mu_s(T))ds<\infty.$$
As is usual, let $\mathcal L_p$ be the submodule of all $T\in\Mcalb$ such that
$$ \int_0^{\infty}\mu_s(T)^p\,ds<\infty.$$

If $\mathcal I$ is a submodule of $\overline{\mathcal M}$ we say
that $\mathcal I$ is {\it geometrically stable} if $\mathcal
I\subset \Mcal+ \mathcal L_{\log}$ and if whenever $h\in \mu(\mathcal
I)$ then $g\in\mu(\mathcal I)$, where
$$ g(t)=\exp\bigg(t^{-1}\int_0^t\log h(s)ds\bigg)\qquad t>0.$$  Geometric
stability is a relatively mild condition.  For example let
$\mathcal X$ be a rearrangement--invariant quasi--Banach function
space on $(0,\infty)$ and suppose $\mathcal I=\{T:\
(\mu_s(T))_{s>0}\in\mathcal X\}\subseteq\Kc+\Lc_{\log}$, where $\Kc\subseteq\Mcal$ is the ideal
of $\tau$--compact operators (see Remark~\ref{rem:obs});
then $\mathcal I$ is
geometrically stable by Proposition 3.2 of \cite{FigielKalton}. A
non--geometrically stable ideal in $B(\HEu)$ is
constructed in \cite{DykemaKalton},
and from this a non--geometrically stable ideal of $\Mcal$ can be constructed.

Suppose $T\in \mathcal L_1\cap \mathcal M$.  Then
the Fuglede--Kadison determinant~\cite{FugledeKadison} of $I+T$ is defined by
$$ \Delta(I+T)=\exp( \tau(\log |I+T|)).$$  Using \cite{Brown} Remark 3.4 we
note that $T\mapsto\log\Delta(I+T)$ is plurisubharmonic on $\Lc_1\cap\Mcal.$
In the Appendix of \cite{Brown} the definition of
$\Delta(I+T)$ is extended to $\mathcal L_{\log}$ and it is shown
that $T\mapsto \Delta(I+T)$ is upper--semicontinuous for the natural
topology of $\mathcal L_{\log}.$  It is not shown explicitly that
$T\mapsto\log\Delta(I+T)$ is plurisubharmonic on $\mathcal L_{\log}$ but this
follows trivially from the results of \cite{Brown}:

\begin{lemma} \label{pluri}  Suppose $S,T\in\Lc_{\log}$.  Then
$$ \log\Delta(I+S)\le \frac{1}{2\pi}\int_0^{2\pi}
\log\Delta(I+S+e^{i\theta}T)\,d\theta.$$\end{lemma}

\begin{proof}  Let $S=H+iK$ and $T=H'+iK'$ be the splitting of
$S,T$ into real and imaginary parts.  Let $R=|H|+|H'|+|K|+|K'|.$
Then $R\in \mathcal L_{\log}$ and $(I+S+zT)(I+R)^{-1}\in I+\mathcal
L_1\cap \mathcal M$ for all $z$.  Using the fact that
$T\mapsto\log\Delta(I+T)$ is plurisubharmonic on $\mathcal L_1\cap \mathcal M$ and
$$ \Delta((I+S+zT)(1+R)^{-1})=\Delta(I+S+zT)(\Delta(I+R))^{-1},$$
it is easy to deduce the Lemma.\end{proof}

Let $g_0(w)=(1-w)$ and
$$g_k(w)=(1-w)\exp(w+\cdots+\frac{w^k}{k})$$ for $k\ge 1.$
If $T\in\Lc_{\log}$ let $k=0$; if $T\in\mathcal M\cap \mathcal L_p$
for some $p>0$, let $k$ be an integer such that $k+1\ge p$.  Then, following \cite{Brown}, there
is a unique $\sigma-$finite measure $\nu=\nu_T$ on $\Cpx\setminus \{0\}$ such that
$$ \log\Delta(g_k(wT))=\int \log|g_k(wz)|d\nu_T(z) \qquad w\in
\Cpx.$$   $\nu_T$ is called the {\it Brown measure} of $T$, and is
independent of the choice of $k$ when many choices are
permissible.  If $T\in \mathcal L_{\log}\cup\bigcup_{p>0}(\mathcal
L_p\cap \mathcal M)$ we shall say that $T$ {\it admits a Brown
measure}.  The measure $\nu_T$ satisfies the following estimates.
If $T\in \mathcal L_{\log}$ and $k=0$ then
\begin{equation}\label{est1} \int_{\Cpx}\log(1+|z|)d\nu_T(z)<\infty
\end{equation}
while if $T\in\mathcal L_p\cap\Mcal$ and $k+1\ge p$, then
\begin{equation}\label{est2} \int_{\Cpx}|z|^pd\nu_T(z)<\infty.
\end{equation}
We refer to \cite[Theorem 3.6]{Brown} and the remark on p.29 of \cite{Brown}.

Of course if $T$ is {\it normal} there is a projection--valued
spectral measure $B\to E_T(B)$ defined for Borel subsets $B$ of the
complex plane and we can define a spectral measure $\nu_T$ by
$$ \nu_T(B)= \tau(E_T(B)).$$  If $T$ also admits a Brown measure
then $\nu_T$ coincides with the Brown measure.

If $T$ either admits a Brown measure or is normal and satisfies $\lim_{t\to\infty}\mu_t(T)=0$,
then for every
$0<r<s<\infty$ we define \begin{equation}\label{Phidef}
\Phi(r,s;T)= \int_{r<|z|\le s} z\,d\nu_T(z).\end{equation}  If $T$
is normal then we can rewrite (\ref{Phidef}) in the form
\begin{equation}\label{Phidef2}
\Phi(r,s;T)= \tau(TE_{|T|}(r,s]).\end{equation}
Note that it is elementary that if $|\alpha|=1$ then
$\Phi(r,s;\alpha T)=\alpha \Phi(r,s;T)$.

\begin{prop}\label{basicprops}
Let $0<r<s<\infty$.
\renewcommand{\labelenumi}{(\arabic{enumi})}
\begin{enumerate}
\item
Suppose $T_1,\ldots,T_N$ are normal with $\lim_{t\to\infty}\mu_t(T_j)=0$ and
$T_1+\cdots+T_N=0$.  Then \begin{equation}\label{qadditive}
|\sum_{j=1}^N \Phi(r,s;T_j)|\le 2N \sum_{j=1}^N(r \tau(
E_{|T_j|}(r,\infty))+s\tau(E_{|T_j|}(s,\infty))).\end{equation}

\item
Suppose $|\alpha|\le 1$ and $T$ is normal with $\lim_{t\to\infty}\mu_t(T)=0$.  Then
\begin{equation}\label{qmult} |\Phi(r,s;\alpha T)-\alpha\Phi(r,s; T)|\le
|\tau(rE_{|T|}(r,\infty)+sE_{|T|}(s,\infty))|.\end{equation}

\item
If $T$ is normal with $\lim_{t\to\infty}\mu_t(T)=0$,  then
\begin{equation}\label{realpart} |\Phi(r,s;\RealPart T)-\RealPart
\Phi(r,s;T)|\le
\tau(rE_{|T|}(r,\infty)+sE_{|T|}(s,\infty))\end{equation} and
\begin{equation}\label{imagpart} |\Phi(r,s;\ImagPart T)-\ImagPart
\Phi(r,s;T)|\le
\tau(rE_{|T|}(r,\infty)+sE_{|T|}(s,\infty))\end{equation}
\end{enumerate}
\end{prop}

\begin{proof} (1) Pick a projection $P\ge E_{|T_j|}(s,\infty)$ for
$1\le j\le N$ and such that $\tau(P)\le
\sum_{j=1}^N\tau(E_{|T_j|}(s,\infty)).$  Then choose $Q\ge P$ with
$Q\ge E_{|T_j|}(r,\infty)$ for $1\le j\le N$ and
$$\tau(Q) \le \sum_{j=1}^N(\tau
(E_{|T_j|}(r,\infty))+\tau(E_{|T_j|}(s,\infty)))\le 2 \sum_{j=1}^N
\tau(E_{|T_j|}(r,\infty)) .$$ Then
$$ \|(Q-E_{|T_j|}(r,\infty))T_j\|\le r,\quad
\|(P-E_{|T_j|}(s,\infty))T_j\|\le s, \qquad 1\le j\le N.$$ Hence
$$ |\tau ((Q-E_{|T_j|}(r,\infty))T_j)|\le r\tau(Q)$$ and
$$|\tau((P-E_{|T_j|}(s,\infty))T_j)| \le s\tau(P).$$
We thus have
\begin{align*}
 |\sum_{j=1}^N \Phi(r,s;T_j)|
&=|\sum_{j=1}^N\tau(T_j(Q-E_{|T_j|}(s,\infty))-T_j(P-E_{|T_j|}(r,\infty)))| \\
&\le N(r\tau(Q)+s\tau(P)).
\end{align*}
Now (\ref{qadditive}) follows.

For (2), we note that
$$ |\Phi(r,s;\alpha T)-\alpha \Phi(r,s;T)|\le
|\alpha|\left(\int_{r<|z|\le|\alpha|^{-1}
r}|z|d\nu_T(z) + \int_{s<|z|\le |\alpha|^{-1}
s}|z|d\nu_T(z)\right).$$  Then (\ref{qmult}) follows
immediately.

Part (3) is similar to (2).  For example we observe for
(\ref{realpart}) that
$$ |\Phi(r,s;\RealPart T)-\RealPart \Phi(r,s;T)|\le \int_{|\RealPart z|\le r< |z|}
|\RealPart z|d\nu_T(z) + \int_{|\RealPart z|\le
s<|z|} |\RealPart z|d\nu_T(z).$$
\end{proof}

\begin{prop}\label{normal}  Let $\mathcal I$ be a submodule of
$\overline{\mathcal M}.$
Suppose $T\in\mathcal I$ is normal
and satisfies $\lim_{t\to\infty}\mu_t(T)=0$.
Then
$T\in[\mathcal I,\mathcal M]$ if and only if there exists a
positive operator $V\in\Ic$ such that
\begin{equation}\label{comm}|\Phi(r,s;T)|\le r\tau(E_V(r,\infty))+s\tau(E_V(s,\infty)) \qquad
0<r<s<\infty.\end{equation}
\end{prop}

\begin{proof}
Assume that (\ref{comm}) holds.
By replacing $V$ with $V+|T|$, if necessary, we may without loss of generality assume
$V\ge|T|$.
Let $h(t)=\mu_t(V).$
Then $h(t)\ge\mu_t(T)$.
If
$0<t<s<\infty$, then from~\eqref{eq:tauE} we have
$$ |\tau(TE_{|T|}(h(s),h(t)])| \le sh(s) +th(t).$$
Now using~\eqref{eq:tauE} again, we get
\begin{align*}
&\big|\tau\big(TE_{|T|}(\mu_s(T),\mu_t(T)]\big)
-\tau\big(TE_{|T|}(h(s),h(t)]\big)\big| \le \\
& \le \int_{\mu_s(T)<|z|\le
h(s)}|z|d\nu_T(z) +\int_{\mu_t(T)<|z|\le h(t)}|z|d\nu_T(z)\\ &\le
sh(s)+th(t).\end{align*}
Hence
$$ \big|\tau\big(TE_{|T|}(\mu_s(T),\mu_t(T)]\big)\big| \le 2sh(s)+2th(t)$$
and we can apply Theorem \ref{thm:IIinf} and (\ref{eq:necsuf}) to
conclude that $T\in [\mathcal I,\mathcal M].$

Conversely, suppose $T$ satisfies (\ref{eq:necsuf}) for some $h$.
Replacing $h$ with
\[
\htil(t)=\frac2t\int_{t/2}^th(s)\dif s,
\]
if necessary, we may without loss of generality
assume $h$ is continuous.
Let $V\in\Ic$  be a positive operator such that $\mu_t(V)=h(t)$.
Given $0<r<s<\infty$, choose $0<v<u$ so that $h(2u)\le r<h(u)$ and
$h(2v)\le s<h(v).$  Then
$$ \big|\tau\big(TE_{|T|}(\mu_{2u}(T),\mu_{2v}(T)]\big)\big|\le 2ur+2vs.$$
Now arguing as above,
\begin{align*} &\big|\tau\big(TE_{|T|}(\mu_{2u}(T),\mu_{2v}(T)]\big)
-\tau\big(TE_{|T|}(r,s]\big)\big| \le \\ & \le \int_{\mu_{2u}(T)<|z|\le
r}|z|d\nu_T(z) +\int_{\mu_{2v}(T)<|z|\le s}|z|d\nu_T(z)\\ &\le
2ur+2vs.\end{align*}
Using Lemma~\ref{lem:infsup}, we have $\tau(E_V(r,\infty))\ge u$ and
$\tau(E_V(s,\infty))\ge v$.
Combining gives
$$\big|\tau\big(TE_{|T|}(r,s]\big)\big|\le 4ur+4vs \le
4r\tau(E_{V}(r,\infty))+4s\tau(E_{V}(s,\infty)).$$
Replacing $V$ by $V\oplus V\oplus V\oplus V$, (cf~ Definition~\ref{defi:oplus}) we have~\eqref{comm}.
\end{proof}

\begin{lemma}\label{subharmonic}  Let $\psi:\Cpx\to\Reals$
be a subharmonic function such that $\psi$ vanishes in a
neighborhood of $0$, is harmonic outside some compact set, and for
a suitable constant $C$, satisfies the estimate  $|\psi(z)|\le
C\log(1+|z|)$ for all $z.$ If $T$ admits a Brown measure, then
define
$$ \Psi(T)=\int_{\Cpx}\psi(z)d\nu_T(z).$$
Suppose $S,T\in \mathcal L_{\log}$ or if $S,T\in\mathcal
L_p\cap\mathcal M$ for some $p>0$.
Then $\Psi(S+e^{i\theta}T)$ is a Borel function of $\theta$ and
\begin{equation}\label{subh}
\Psi(S) \le \frac{1}{2\pi}\int_0^{2\pi}
\Psi(S+e^{i\theta}T)d\theta.\end{equation}\end{lemma}

\begin{proof} By an easy approximation argument it will suffice to
consider the case when $\psi$ is $C^2.$  In this case for any
choice of $k\ge 0$ we have the formula (\cite{Brown} Proposition
2.2)
$$ \psi(z)=\int_{\Cpx}\log|g_k(w^{-1}z)|\nabla^2\psi(w)
d\lambda(w) \qquad z\in\Cpx$$ where $\lambda$ denotes area
measure.  Hence if $T$ admits a Brown measure and $k$ is suitably
chosen, \begin{equation}\label{fubini} \Psi(T)=\int_{\Cpx}(
\int_{\Cpx}\log|g_k(w^{-1}z)|
\nabla^2\psi(w)d\lambda(w))d\nu_T(z).\end{equation} Now it can be
checked that the function $|\log g_k(w^{-1}z)|\nabla^2\psi(w)$ is
integrable for the product measure $\lambda\times \nu_T.$
Indeed, let us first
consider the case when $T\in\mathcal L_p\cap\mathcal M$, with
$k+1\ge p$.
Estimates on the growth of $\log|g_k(w)|$ (cf. p.~11 of~\cite{Brown}) give
$$ \int_0^{2\pi} \big|\log |g_k(r^{-1}e^{-i\theta}z)|\big|d\theta \le
C\min (|z|^{k+1}|r|^{-k-1}, |z|^{k+\epsilon}|r|^{-k-\epsilon})$$
for suitable $C$ and $\epsilon>0$.
Since
$\nabla^2\psi$ has compact support contained in some annulus away
from the origin we need only observe that
$$ \int \min(|z|^{k+1},|z|^{k+\epsilon})d\nu_T(z) <\infty$$ which
follows from (\ref{est2}).
On the other hand, if $T\in \mathcal L_{\log}$ and thus $k=0$,
we use the estimate
$$ \int_0^{2\pi}\big|\log|g_k(r^{-1}e^{-i\theta}z)|\big|d\theta\le
C\log(1+|z|)$$ and (\ref{est1}). It follows we can use Fubini's
theorem to rewrite (\ref{fubini}) in the form
\begin{align*}
\Psi(T)&=\int_{\Cpx}( \int_{\Cpx}\log|g_k(w^{-1}z)|
d\nu_T(z))\nabla^2\psi(w)d\lambda(w))\\
&= \int_{\Cpx}
\log\Delta(g_k(w^{-1}T))\nabla^2\psi(w)d\lambda(w).\end{align*}
Now the result follows easily from the upper semicontinuity of
$\log\Delta$ and Lemma \ref{pluri}.\end{proof}

\begin{prop}\label{normalspectrum}  Let $\mathcal I$ be a
geometrically stable submodule of $\overline{\mathcal M}.$
If $T\in\mathcal I$ admits a Brown measure, then there is a normal
operator $S\in\mathcal I$ with $\nu_S=\nu_T$.
Furthermore, $S$ admits a Brown measure.
\end{prop}

\begin{proof}  It will suffice to show the existence of a positive
operator $V\in \mathcal I$ so that
$$ \nu_T(|z|>r) \le \nu_V(r,\infty) \qquad 0<r<\infty.$$
Let $H=\RealPart T,\ K=\ImagPart T$ and then set $P=|H|+|K|.$
Since $\mathcal I$ is geometrically stable there exists a positive
$V\in\mathcal I$ with
$$ \frac1t\int_0^t\log\mu_s(P)ds \le \log\mu_t(V) \qquad 0<t<\infty.$$
Therefore, $\mu_t(P)\le\mu_t(V)$ and $\nu_P(r,\infty)\le\nu_V(r,\infty)$ for all $0<r<\infty$.

Suppose for contradiction that for some $0<r<\infty$ we have
$t=\nu_T(|z|>r)>\nu_V(r,\infty)$.  Choose
$r_0<r$ so that $\nu_P[r_0,\infty)\ge t\ge \nu_P(r_0,\infty).$
Let $\psi(z)=\log_+\frac{|z|}{r_0}$ and define $\Psi$ as in Lemma
\ref{subharmonic}.  Then
$$ \Psi(T) \le \frac{1}{2\pi}\int_0^{2\pi}
\Psi(T+e^{i\theta}T^*)d\theta.$$ Now $T+e^{i\theta}T^*=
2e^{i\theta/2}(H\cos\frac{\theta}2+K\sin\frac{\theta}2).$  Hence
$|T+e^{i\theta}T^*|\le2(|H|+|K|)=P$ and it follows that
$\Psi(T+e^{i\theta}T^*)\le \Psi(P)$ for $0\le\theta\le 2\pi.$
\begin{equation*}
t\log\frac{r}{r_0}\le
\int\log_+\frac{|z|}{r_0}d\nu_T(z)=\Psi(T)
\le \Psi(P)
= \int_0^t\log_+\frac{\mu_s(P)}{r_0}ds
\le t\log \frac{\mu_t(V)}{r_0}.
\end{equation*}
Thus $\mu_t(V)\ge r$
and hence $\nu_V(r,\infty)\ge t$ contrary to assumption.

The inequalities~\eqref{est1} and~\eqref{est2} imply that $S$ admits a Brown measure.
\end{proof}

Before proving our main result it will be convenient to introduce
some notation.  Let $\mathcal I$ be any submodule of
$\overline{\mathcal M}$ not containing $\Mcal$.
Hence $\lim_{t\to\infty}\mu_t(T)=0$ for every $T\in\Ic$.
Let $F(r,s)$ be a function of two
variables defined for $0<r<s<\infty.$  We write $F\in\mathcal
F(\mathcal I)$ if there exists a positive operator $V\in\mathcal
I$ such that
$$ |F(r,s)|\le r\tau (E_V(r,\infty))+s\tau (E_V(s,\infty))\qquad 0<r<s<\infty.$$
We write $F\in\mathcal G(\mathcal I)$ if there if there is a
positive operator $V\in\mathcal I$ such that
$$ |F(r,s)|\le \int_{(0,\infty)}\big(r\log_+\frac xr
+s\log_+\frac xs\big)\,d\nu_V(x) \qquad 0<r<s<\infty.$$ Both
$\mathcal F(\mathcal I)$ and $\mathcal G(\mathcal I)$ are easily
seen to be vector spaces.  Also note that $\mathcal F(\mathcal
I)\subset \mathcal G(\mathcal I)$ (replace $V$ by $eV$.)
Proposition \ref{normal} states that if $T$ is normal then $T\in
[\mathcal I,\mathcal M]$ if and  only if $\Phi(r,s;T)\in\mathcal
F(\mathcal I).$  We improve this for geometrically stable
submodules.

\begin{prop}\label{hermitian}  Suppose $\mathcal I$ is a
geometrically stable submodule of $\overline{\mathcal M}$ with $\Mcal\not\subseteq\Ic$.  If
$T\in\mathcal I$ is normal,
then $T\in [\mathcal I,\mathcal M]$
if and only if $\Phi(r,s;T)\in\mathcal G(\mathcal I).$\end{prop}

\begin{proof} One direction is trivial from Proposition
\ref{normal}.  For the other direction suppose
$\Phi(r,s;T)\in\mathcal G(\mathcal I).$  Choose $V$ a positive
operator in $\mathcal I$ so that
\begin{equation}\label{eq:Phirs}
 |\Phi(r,s;T)|\le \int_{(0,\infty)}
\big(r\log_+\frac xr+s\log_+\frac xs\big)\,d\nu_V(x) \qquad
0<r<s<\infty.
\end{equation}
We can assume $V\ge |T|.$ Let $h(t)=\mu_t(V)$ and
let $$g(t)=\exp(\frac1t\int_0^t\log h(s)ds)\qquad 0<t<\infty.$$
Suppose $0<t<s<\infty.$ Then similarly to in the proof of Proposition~\ref{normal}, we get
$$ \big|\tau\big(TE_{|T|}(\mu_s(T),\mu_t(T)]\big)\big|\le
\big|\tau\big(TE_{|T|}(\mu_s(V),\mu_t(V)]\big)\big| + sh(s) +th(t).$$
Now from~\eqref{eq:Phirs} we get
\begin{align*} \big|\tau\big(TE_{|T|}(\mu_s(V),\mu_t(V)]\big)\big|&\le \int_{(0,\infty)}
\big(\mu_s(V)\log_+\frac x{\mu_s(V)}+
\mu_t(V)\log_+\frac x{\mu_t(V)}\big)d\nu_V(x)\\
&= \int_0^s h(s)\log\frac{h(u)}{h(s)}du +\int_0^t
h(t)\log\frac{h(u)}{h(t)}du\\
&= sh(s)\log \frac{g(s)}{h(s)}+th(t)\log\frac{g(t)}{h(t)}\\
&\le sg(s)+tg(t).\end{align*} Combining, we see that
$$ \big|\tau\big(TE_{|T|}(\mu_s(T),\mu_t(T)]\big)\big|\le
s(h(s)+g(s))+t(h(t)+g(t))$$ and so by Theorem \ref{thm:IIinf},
$T\in\mathcal [I,\mathcal M].$\end{proof}

\begin{thm}\label{real}  Suppose $\mathcal I$ is a submodule of
$\overline{\mathcal M}$ with $\Mcal\not\subseteq\Ic$ and $T\in\mathcal I$ admits a Brown
measure.   Then
$$ \RealPart\Phi(r,s;T)-\Phi(r,s;\RealPart T)\in \mathcal G(\mathcal
I),\quad \ImagPart \Phi(r,s;T)-\Phi(r,s;\ImagPart T)\in\mathcal
G(\mathcal I).$$
\end{thm}

\begin{proof}  Let $H=\RealPart T$ and $K=\ImagPart T.$ We need only prove the statement concerning the real
part, since the other half follows by considering $iT.$ We also
note that if $s\le 2r$ we have $|\Phi(r,s,T)|\le 2r\nu_T(|z|>r)$
and $|\Phi(r,s;H)|\le 2r\nu_{|H|}(r,\infty).$  By Proposition~\ref{normalspectrum},
this implies an estimate
$$|\RealPart\Phi(r,s;T)-\Phi(r,s;H)|\le 2r\nu_V(r,\infty) \qquad 0<r<s\le 2r<\infty$$ for a
suitable positive operator $V\in\mathcal I.$  This means we need
only consider estimates when $s>2r.$

We first fix a smooth bump function $b:\Reals\to\Reals$
such that $\supp b\subset (0,1/2)$, $b\ge 0$, $\int b(x)dx=1.$ Let
$\beta(t)=2|b(t)|+|b'(t)|.$

Now suppose $0<r<s<\infty$, with $s> 2r$.  We define
$$ \varphi_{r,s}(\tau) = \int_{-\infty}^{\tau} b(t-\log
r)- b(t-\log s)dt.$$  Notice that the two terms in the integrand
are never simultaneously positive (since $\log 2>\frac12$), and
$\varphi_{r,s}$ is a bump function which satisfies
$\varphi_{r,s}(\tau)=0$ if $\tau<\log r$ or $\tau>\frac12+\log s$,
while $\varphi_{r,s}(\tau)=1$ if $\frac12+\log r \le \tau\le \log
s$ and $0\le\varphi_{r,s}(\tau)\le1$ for all $\tau$.

Then let $\rho_{r,s}$ be defined to be the function such that
$\rho_{r,s}(\tau)=0$ if $\tau<\log r$ and
$$\rho''_{r,s}(\tau)=
e^{\tau}(2|\varphi_{r,s}'(\tau)|+|\varphi_{r,s}''(\tau)|).$$ In
fact, this implies that $$\rho_{r,s}''(\tau)=
e^{\tau}(\beta(\tau-\log r)+\beta(\tau-\log s))$$ and then
$$ \rho_{r,s}'(\tau)= \int_{-\infty}^{\tau}e^t (\beta(t-\log r)+\beta(t-\log
s))dt$$
and
$$ \rho_{r,s}(\tau)=\int_{-\infty}^{\tau}(\tau-t)e^t
(\beta(t-\log r)+\beta(t-\log s))\,dt.$$
Thus, if we set
$$ C_0=\int_{-\infty}^{\infty} e^t \beta(t)dt,$$
then
$$ \rho_{r,s}'(\tau)\le
C_0(r\chi_{(\tau
>\log r)}
+s\chi_{(\tau>\log s)})$$ and so
\begin{equation}\label{upperpsi} 0\le\rho_{r,s}(\tau) \le C_0(r (\tau-\log
r)_++s(\tau-\log s)_+).\end{equation}

Now we use the argument of Lemma 2.6 of \cite{Kreine}.   We define
$$ \psi_{r,s}(z)= \rho_{r,s}(\log |z|)-x\varphi_{r,s}(\log |z|) \qquad
z=x+iy\neq 0$$ and $\psi(0)=0$.  Then if $z\neq 0,$
$$\nabla^2\rho_{r,s}(\log|z|)= |z|^{-2}\rho''_{r,s}(\log|z|).$$ Similarly
$$ \nabla^2(x\varphi_{r,s}(\log|z|)=
\frac{x}{|z|^2}(2\varphi'_{r,s}(\log|z|)+\varphi''_{r,s}(\log|z|)).$$
Thus by construction, $|\nabla^2(x\varphi_{r,s}(\log|z|)|\le
\nabla^2 (\rho_{r,s}(\log|z|)$ and so $\psi_{r,s}$ is subharmonic.
Note that $\psi_{r,s}$ also vanishes on a neighborhood of $0$ and
is harmonic outside a compact set. We note the estimates (from
(\ref{upperpsi}))
\begin{equation}\label{upperpsi2} 0\le\rho_{r,s}(\log |z|)\le
C_0\left(r\log_+\frac{|z|}{r}+s\log_+\frac{|z|}{s}\right)\end{equation}
and
\begin{equation} \label{hupper} 0\le \psi_{r,s}(z) \le  C_0\left(r\log_+
\frac{|z|}{r}+s\log_+\frac{|z|}{s}\right) \qquad |z|\ge 2s.
\end{equation} Note of course that $C_0$ is independent of $r,s.$

If $A$ admits a  Brown measure or is normal with $\lim_{t\to\infty}\mu_t(A)=0$, let us define
\begin{align*}
\tilde\Phi(r,s;A)&=\int_{\Cpx}(\RealPart z)\varphi_{r,s}(\log|z|)d\nu_A(z) \\[1ex]
\Omega(r,s;A)&=\int_{\Cpx} \rho_{r,s}(\log |z|)d\nu_A(z) \\[1ex]
\Psi(r,s;A)&=\int_{\Cpx}\psi_{r,s}(z)d\nu_A(z).
\end{align*}
Thus $\Psi(r,s;A)=\Omega(r,s;A)-\tilde \Phi(r,s;A)$ and
$\Psi(r,s;-A)=\Omega(r,s,A)+\tilde\Phi(r,s;A)$.    We can apply
Lemma \ref{subharmonic} to
$\Psi(r,s;\cdot)$, giving
\begin{align*} \Psi(r,s;T) &\le \frac{1}{2\pi}\int_0^{2\pi}
\Psi(r,s;T+e^{i\theta}T^*)d\theta \\
\Psi(r,s;-T) &\le \frac{1}{2\pi}\int_0^{2\pi}
\Psi(r,s;-T-e^{i\theta}T^*)d\theta.
\end{align*}
Note that
$\theta\to \tilde \Phi(r,s;T+e^{i\theta}T^*)$ is a Borel function
by using Lemma \ref{subharmonic} and the equation
$$
\tilde\Phi(r,s;A)=\frac12(\Psi(r,s;A)-\Psi(r,s;-A)).$$
We have
\begin{equation}\label{111}
\left|\frac{1}{2\pi}\int_0^{2\pi}\tilde\Phi(r,s;T+e^{i\theta}T^*)d\theta-
\tilde\Phi(r,s;T)\right| \le \frac{1}{2\pi} \int_0^{2\pi}
\Omega(r,s;T+e^{i\theta}T^*)d\theta.\end{equation}

We first estimate the right--hand side of (\ref{111}).  Note that
$\Omega(r,s;T+e^{i\theta}T^*)= \Omega(r,s;W_{\theta})$ where
$W_{\theta}=2(H\cos\frac{\theta}2+K\sin\frac{\theta}2)$ is hermitian
and hence from~\eqref{upperpsi2},
$$ \Omega(r,s,W_{\theta})\le
C_0\left(r\int_0^{\infty}\log_+\frac{\mu_t(|W_{\theta}|)}{r}dt +
s\int_0^{\infty}\log_+\frac{\mu_t(|W_{\theta}|)}{s}dt \right).$$
Hence for all $0\le \theta\le 2\pi,$
$$ \Omega(r,s,W_{\theta})\le
C_0\left(r\int_0^{\infty}\log_+\frac{\mu_t(P)}{r}dt +
s\int_0^{\infty}\log_+\frac{\mu_t(P)}{s}dt \right),$$ where
$P=2(|H|+|K|).$  Thus the right--hand side of (\ref{111}) is estimated
by
$$C_0\left(r\int_0^{\infty}\log_+\frac{\mu_t(P)}{r}dt +
s\int_0^{\infty}\log_+\frac{\mu_t(P)}{s}dt \right).$$ In other
words the right--hand side of (\ref{111}) belongs to $\mathcal
G(\mathcal I),$ and hence so does the left--hand side.

Now we turn to the left--hand side of (\ref{111}).    We note that
$$ |\tilde \Phi(r,s;T)-\RealPart\Phi(r,s;T)|\le
\int_{r<|z|<2r}|z|d\nu_T(z)+\int_{s<|z|<2s}|z|d\nu_T(z).$$  Hence
$$ |\tilde \Phi(r,s;T)-\RealPart\Phi(r,s;T)|\le
2r\nu_T(|z|>r)+2s\nu_T(|z|>s).$$  By Proposition
\ref{normalspectrum} this implies that $\tilde
\Phi(r,s;T)-\RealPart \Phi(r,s;T)\in\mathcal F(\mathcal I).$

By the same argument we also have
$$ \sup_{0\le \theta\le
2\pi}|\tilde\Phi(r,s;T+e^{i\theta}T^*)-\RealPart
\Phi(r,s;T+e^{i\theta}T^*)|\in\mathcal F(\mathcal I).$$

Now, by using parts (1) and (3) of Proposition \ref{basicprops}, we easily obtain that
$$ \sup_{0\le \theta\le 2\pi} |\RealPart\Phi(r,s;T+e^{i\theta}T^*)-
(1+\cos\theta)\Phi(r,s;H)+\sin\theta\,\Phi(r,s;K)|\in \mathcal
F(\mathcal I).$$   So on integration we find that
$$
\frac{1}{2\pi}\int_0^{2\pi}\tilde\Phi(r,s;T+e^{i\theta}T^*)d\theta-
\Phi(r,s;H) \in \mathcal F(\mathcal I).$$  It follows that the
left--hand side of (\ref{111}) differs from
$|\RealPart\Phi(r,s;T)-\Phi(r,s;H)|$  by a function in class
$\mathcal F(\mathcal I).$  Combining we obtain:
$$ \RealPart\Phi(r,s;T)-\Phi(r,s;H)\in \mathcal G(\mathcal
I).$$\end{proof}

\begin{thm}\label{brown}  Let $\mathcal I$ be a geometrically
stable submodule of $\overline{\mathcal M}.$  Let $T\in\mathcal I$
admit a Brown measure.  Then $T\in[\mathcal I,\mathcal M]$ if and
only if there is a positive operator $V\in \mathcal I$ with
\begin{equation}\label{eq:BrownEv}
\bigg|\int_{r<|z|\le s}z\,d\nu_T(z)\bigg|\le
r\tau(E_V(r,\infty))+s\tau(E_V(s,\infty)) \qquad 0<r,s<\infty.
\end{equation}
\end{thm}

\begin{proof}
First suppose $\Mcal\not\subseteq\Ic$.
Let $H=\frac12(T+T^*)$ and $K=\frac{1}{2i}(T-T^*)$.
Note that
$T\in\mathcal [I,M]$ if and only if $H,K\in[\mathcal I,\mathcal
M].$
Then by Theorem \ref{real} we have $\Phi(r,s;T)\in \mathcal
G(\mathcal I)$ if and only if $\Phi(r,s;H),\Phi(r,s;K)\in\mathcal
G(\mathcal I).$  By Proposition \ref{hermitian} this implies that
$\Phi(r,s;T)\in \mathcal G(\mathcal I)$ if and only if $T\in[\mathcal
I,\mathcal M].$

Let $S\in\mathcal I$ be a normal operator with $\nu_S=\nu_T$ as
given by Proposition \ref{normalspectrum}.
Then the same reasoning as above applies to $S$,
yielding $S\in[\Ic,\Mcal]$ if and only if $\Phi(r,s;S)\in\Gc(\Ic)$.
By Proposition~\ref{normal}, $S\in[\Ic,\Mcal]$ if and only if $\Phi(r,s;S)\in\Fc(\Ic)$.
But $\Phi(r,s;T)=\Phi(r,s;S)$, so $T\in[\Ic,\Mcal]$ if and only if
$\Phi(r,s;T)\in\mathcal F(\mathcal I)$.

Now suppose $\Mcal\subseteq\Ic$.
If $T\in[\Ic,\Mcal]$, then by Proposition~\ref{prop:I0}, $T\in[\Ic_0,\Mcal]$, so
by the case just proved there is a positive operator $V\in\Ic_0$ making~\eqref{eq:BrownEv} hold.
On the other hand, suppose $V\in\Ic$ is a positive operator making~\eqref{eq:BrownEv} hold.
Let $S\in\mathcal I$ be a normal operator with $\nu_S=\nu_T$ as
given by Proposition \ref{normalspectrum}.
Then
\[
|\Phi(r,s;S)|=|\Phi(r,s;T)|\le
r\tau(E_V(r,\infty))+s\tau(E_V(s,\infty)) \qquad 0<r,s<\infty.
\]
Hence, by Proposition~\ref{normal}, $S\in[\Ic,\Mcal]$.
Invoking Propositions~\ref{prop:I0} and~\ref{normal} again, we find a positive operator
$V'\in\Ic_0$ such that 
\[
|\Phi(r,s;S)|\le
r\tau(E_{V'}(r,\infty))+s\tau(E_{V'}(s,\infty)) \qquad 0<r,s<\infty.
\]
But then, since $\Mcal\not\subseteq\Ic_0$, we get $T\in[\Ic_0,\Mcal]$ by the case proved above.
\end{proof}

Let us say that $T$ is {\it approximately nilpotent} if $T$
admits a Brown measure with $\nu_T=0.$  This is equivalent to the
statement that $\Delta(g_k(wT))=1$ for all $w\in\Cpx.$

\begin{cor}  If $\mathcal I$ is a geometrically stable submodule
of $\overline{\mathcal M}$ then every approximately nilpotent
$T\in\mathcal I$ belongs to $[\mathcal I,\mathcal M].$\end{cor}

We remark that in the case of an ideal $\Ic$ of $B(\HEu)$ the
relationship between the subspace $[\Ic,B(\HEu)]$
and the growth of the characteristic determinant is discussed
further in \cite{kaltonpalermo}, and it is possible that some
analogous results can be obtained here for the Fuglede--Kadison
determinant.

\bibliographystyle{plain}

\end{document}